\documentclass[10pt]{article}
\usepackage{amsmath}
\usepackage{amssymb}
\usepackage{amsthm}
\usepackage{amsfonts}
\usepackage[utf8]{inputenc}
\usepackage{tikz,pgfplots}

\usepackage{fullpage}
\usepackage{tabularx}
\usepackage{multirow}
\usepackage{rotating}
\usepackage{enumerate}

 \usepackage[usenames,dvipsnames]{pstricks}
 \usepackage{epsfig}
 \usepackage{pst-grad} 
 \usepackage{pst-plot} 
 \usepackage[space]{grffile} 
 \usepackage{etoolbox} 
 \makeatletter 
 \patchcmd\Gread@eps{\@inputcheck#1 }{\@inputcheck"#1"\relax}{}{}
 \makeatother

\usepackage{graphicx}  
\usepackage{caption,subcaption}

\usepackage[software,hardware]{mymacros}
\newcommand{\smb}{\left[\begin{smallmatrix}}
\newcommand{\sme}{\end{smallmatrix}\right]}

\newcommand{\beq}{\begin{equation}}
\newcommand{\eeq}{\end{equation}}

\usepackage{mathdots}
\usepackage[colorlinks=true]{hyperref}
\usepackage{algorithm}
\usepackage{algpseudocode} 

\theoremstyle{definition}

\title{Numerical Linear Algebra in Data Assimilation} %
\author{Melina A. Freitag \thanks{Institut f\"{u}r Mathematik, Universit\"{a}t Potsdam, 
Campus Golm, Haus 9, Karl-Liebknecht-Str. 24-25, D-14476 Potsdam, {\tt
  melina.freitag@uni-potsdam.de}}}

\begin{document}
\maketitle

\begin{abstract}
Data assimilation is a method that combines observations (that is, real world data) of a state of a system with model output for that system in order to improve the estimate of the state of the system and thereby the model output. The model is usually represented by a discretised partial differential equation. The data assimilation problem can be formulated as a large scale Bayesian inverse problem. Based on this interpretation we will derive the most important variational and sequential data assimilation approaches, in particular three-dimensional and four-dimensional variational data assimilation (3D-Var and 4D-Var) and the Kalman filter. We will then consider more advanced methods which are extensions of the Kalman filter and variational data assimilation and pay particular attention to their advantages and disadvantages. The data assimilation problem usually results in a very large optimisation problem and/or a very large linear system to solve (due to inclusion of time and space dimensions). Therefore, the second part of this article aims to review advances and challenges, in particular from the numerical linear algebra perspective, within the various data assimilation approaches.

\end{abstract}

\section{Motivation}
\label{sec:intro}

Integrating large data sets into sophisticated computational models is one of the big challenges in mathematical sciences of the 21st century.
When the computational model arises from a dynamical system, and time dependent observational data of that system are available, then the process of combining the model and the data to obtain a more informed system is called \emph{data assimilation}.

Data assimilation research has been mainly driven by practitioners, initially in the field of numerical weather prediction and ocean modelling \cite{ghil1991,Ghil1989,Courtier1994,Courtier1993,
Sasaki1970,Sasaki1958,Talagrand1987,daley1993,nichols2003,Lawless2013,navon2009}, but nowadays has many more applications in geosciences \cite{carrassi2018,vanLeeuwen2010,fletcher2017,park2013}, ecology \cite{niu2014,luo2011}, biology \cite{lawson1995,rhodes2009}, chemistry \cite{brown1986,elbern1997}, mechanical engineering \cite{altaf2013,corigliano2004}, medicine \cite{kostelich2011,Engbert2020}, image processing \cite{bereziat2011,butala2009}, as well as human and social sciences \cite{schiff2009,schwartz2014}, see also \cite{Asch2016} and references therein, with the potential for further utilisation in data science
and machine learning. In particular, as it becomes easier to make large numbers of relatively accurate observations of a system (we explain later what we mean by a system), a major challenge is how best to use this information to update and refine the model of that system. Only in recent years a more mathematical approach to the theory of data assimilation has been developed, see for example \cite{Stu2010,kaipio2006statistical,wikle2007,LawStuZyg2015,FrePot2013,ReiCot2013}. Good introductions to Bayesian approaches for inverse problems from the viewpoint of numerical linear algebra are the book \cite{calv2007} and the reference \cite{calv2018}.

A model of a system consists of a number of mathematical equations, often (stochastic) partial differential equations that describe the interactions between several variables through certain processes from physics, biology, chemistry or other applications. Most of the time those equations simplify the dynamics of the system by excluding processes that are considered less important, happen at different scales, or are simply not easy to model. Moreover, parameters in the system may only be known approximately, and the computational model (via discretisation of differential equations)  results in another error introduced within the process. Even if we did have perfect knowledge of a system, often initial conditions or boundary conditions are not known to high accuracy. 

In addition to the model, measurements of the system variables, perhaps indirect, are available at different locations in time and space. 

Both the model and the observations obtained through measurements have errors. The data assimilation problem is therefore an inverse problem that uses incomplete and erroneous data and knowledge about a model which is also imperfect, in order to find the best possible estimate of the state of a system (or the best possible approximation of an unknown parameter, in which case we speak of parameter estimation). Using this state estimate one can then use the model to make predictions about the future states of the system, and, as more observations become available, update the state estimate using cycled data assimilation schemes.

Many algorithms for data assimilation have been developed over the past century, we introduce the main ones in the next section. However many challenges remain, in particular as more and more data become available and larger, more complex and higher dimensional problems need to be solved.

This review article will not address all problems arising within data assimilation, but will focus on challenges in data assimilation for numerical linear algebra.

The article is structured as follows: The most common methods for data assimilation are introduced in Section \ref{sec:basic}: Based on Bayes' theorem we derive variational and sequential data assimilation techniques. Section \ref{sec:lin} focuses on the solution to the optimisation problem and the linear system arising within variational data assimilation. Approximations, in particular low-rank approximations, and other variants of the Kalman filter are considered in Section \ref{sec:kfred}. Section \ref{sec:modred} reviews several dimension reduction approaches for variational and sequential data assimilation methods, and in Section \ref{sec:baytic} we briefly consider various other aspects, such as the connection between data assimilation, Bayesian inference and Tikhonov regularisation. Finally the last two sections give a short survey on data assimilation software and a conclusion. By no means we consider this article a complete survey on data assimilation techniques, we rather focus on a selected methods and approaches where linear algebra plays an important role.

\section{Basic methods and algorithms for state estimation}
\label{sec:basic}

The goal of data assimilation is to incorporate measured observations into a model of a dynamical system in order to produce estimates of the current system state (and future system states) which are as accurate as possible. In that sense data assimilation can be defined as an approximation of a true state of a physical system at a given time, by combining time-distributed observations with the dynamical system model in some optimal way.

One can view the data assimilation problem as a Bayesian inference problem. Let $x\in\mathbb{R}^n$ be a model state that we would like to estimate. In Bayesian statistics, we model $x$ as a realisation of a random variable (here a random vector), $X:\Omega\rightarrow\mathbb{R}^n$. If $\Theta:\Omega\rightarrow\mathbb{R}^p$ is another random variable with mean zero modelling the observational noise, then we can also model the observed variable $Y\in\mathbb{R}^p$ as a random variable, defined by
\[
Y = h(X)+\Theta,
\]  
where $h:\mathbb{R}^n\rightarrow\mathbb{R}^p$ is a (in general nonlinear) continuous map which models the transformation of the system space to the observation space. $X$ and $\Theta$ are assumed to be independent. We would like to infer information about states $x$ given realisations $y$ of $Y$, which is a Bayesian inverse problem \cite{bardsley2018,Stu2010,dashti2016,allmaras2013,calv2007,calv2018,apte2008}.

If we assume that $\Theta\in\mathbb{R}^p$ is a random variable with probability density $\pi$, and $y$ as well as $x$ realisations of the random variable $Y$ and $X$, respectively, then the probability of $y$ given $x$ is given by $\pi_{Y|X}(y|x) = \pi(y-h(x))$. This is often referred to as the data likelihood. Further let $\pi_X(x)$ be the probability density function of $X$, which describes our prior believes about the distribution of $X$. Then, by Bayes' formula, $\pi_{X|Y}(x|y)$, the \emph{posterior conditional probability density} function of $x$ given the observations $y$ is given by
\beq
\label{eq:bayes}
\pi_{X|Y}(x|y) = \frac{\pi_{Y|X}(y|x)\pi_X(x)}{\pi_Y(y)}\propto\pi_{Y|X}(y|x)\pi_X(x),
\eeq
where $\pi_Y(y) = \int_{\mathbb{R}^n}\pi_{Y|X}(y|x)\pi_X(x)dx$ is a normalisation constant depending only on $y$ (see, for example \cite{Stu2010,LawStuZyg2015,reich2015}). In general it is hard to obtain the entire probability density $\pi_X(x|y)$, in particular in higher dimensions. However, if we make some assumptions about the probability density functions of the prior and the likelihood, the problem of finding the posterior density can be simplified, which leads to data assimilation algorithms in the traditional sense, which we discuss in the next sections. We will distinguish between variational and sequential data assimilation methods.

\subsection{Variational data assimilation}
\label{sec:var}

If the prior density in the Bayesian inference problem (\ref{eq:bayes}) is Gaussian, that is $X\sim\mathcal{N}(x^B,B)$ with mean $x^B\in\mathbb{R}^n$ (often called the background vector) and positive definite background error covariance matrix $B\in\mathbb{R}^{n\times n}$, then we have  
\[
\pi_X(x) = \frac{1}{\sqrt{(2\pi)^n \det{B}}}\exp\left(-\frac{1}{2}(x-x^B)^TB^{-1}(x-x^B)\right).
\]
Here $x\in\mathbb{R}^n$ is the state vector we would like to estimate. If, in addition, the observation error is also Gaussian, that is $\Theta\sim\mathcal{N}(0,R)$ (or $Y\sim\mathcal{N}(h(X),R)$) with positive definite error covariance matrix $R\in\mathbb{R}^{p\times p}$, then the likelihood function can be written as 
\[
\pi_{Y|X}(y|x) = \frac{1}{\sqrt{(2\pi)^p \det{R}}}\exp\left(-\frac{1}{2}(y-h(x))^T R^{-1}(y-h(x))\right).
\]
and hence 
\[
\pi_X(x|y)\propto\exp\left(-\frac{1}{2}\|y-h(x)\|^2_{R^{-1}} -\frac{1}{2}\|x-x^B\|^2_{B^{-1}}\right),
\]
where $\|z\|^2_{R^{-1}} := z^T R^{-1}z$ is a weighted norm, the so-called Mahalanobis distance to zero. The maximum a posteriori (MAP) estimator for $x\in\mathbb{R}^n$ is then given by
\beq
\label{eq:3dvar}
\underset{x\in\mathbb{R}^n}{\operatorname{argmin }} J(x), \quad\text{where}\quad J(x) =  \left(\frac{1}{2}\|y-h(x)\|^2_{R^{-1}} +\frac{1}{2}\|x-x^B\|_{B^{-1}}^2\right),
\eeq
which is a weighted nonlinear least squares problem. The minimisation problem in (\ref{eq:3dvar}) is what is known as the three dimensional variational data assimilation problem (\emph{3D-Var}). If $h:\mathbb{R}^n\rightarrow\mathbb{R}^p$ is a linear operator, represented by a matrix $H\in\mathbb{R}^{p\times n}$, the solution to (\ref{eq:3dvar}) can be computed immediately using the (sometimes called Kalman gain) matrix $K$ given by
\[
K = BH^T(HBH^T+R)^{-1},
\]
or $K = (B^{-1}+H^TR^{-1}H)^{-1} H^T R^{-1}$, using the Sherman-Morrison-Woodbury formula \cite{Golub2012}, and the minimum in \eqref{eq:3dvar} is then given by
\[
x^* = x^B+K(y-H(x^B)).
\]
This direct solve of the optimisation problem \eqref{eq:3dvar} is sometimes referred to as optimal interpolation \cite{daley1993}. 

We note that the minimisation problem in (\ref{eq:3dvar}) is a generalised form of Tikhonov regularisation \cite{FrePot2013,MuellSilt2012,bardsley2018,hansen2010,neumaier1998,Johnson2005}, the most commonly used form of regularisation for inverse problem. In standard Tikhonov regularisation problems, the operator $R=I$, where $I$ is the identity matrix, $x^B=0$, $B = \frac{1}{\lambda} I$, where $\lambda>0$ is a regularisation parameter, and often $h$ is linear. The generalised form (\ref{eq:3dvar}) is a nonlinear Tikhonov regularisation in a  weighted inner product space (see, \cite{FrePot2013,calv2018,vogel2002}). In recent years the dynamic aspect of inverse problems (naturally leading to a variety of data assimilation problems) has become of interest in the inverse problems community, we especially refer to \cite{Schuster_2018,judd2008forecasting} and references therein.

In practice the matrices involved in the computation of the Kalman gain are very large and direct inversion is not feasible, often it is not even possible to store the full matrices. In addition, if the operator $h$ is nonlinear, then an iterative optimisation procedure is required in order to solve the minimisation problem (\ref{eq:3dvar}), we will discuss more details about this optimisation in Section \ref{sec:lin}.

In variational data assimilation one uses a descent algorithm in the direction of the gradient and an adjoint approach for the computation of the gradient in order to solve the minimisation problem. 

The problem (\ref{eq:3dvar}) is a stationary one, there is no time-dependence. Instead of a state $x$ fixed in time, let us consider a set of state estimates $x=[x^T_0,\ldots, x^T_N]^T$, where $x_i\in\mathbb{R}^n$ refers to a state at time $t_i$. Define a new observation operator $\mathcal{H}:\mathbb{R}^{n(N+1)}\rightarrow\mathbb{R}^{p(N+1)}$, with $\mathcal{H}:[x^T_0,\ldots, x^T_N]^T\rightarrow [\mathcal{H}_0(x_0)^T,\ldots, \mathcal{H}_N(x_N)^T]^T$ and let $y=[y^T_0,\ldots, y^T_N]^T$, where $y_i\in\mathbb{R}^p$ denotes a set of observations at time $t_i$, $i = 0,\ldots,N$. Furthermore, let $R\in\mathbb{R}^{p(N+1)\times p(N+1)}$ now be a block diagonal matrix with $R_i\in\mathbb{R}^{p\times p}$, $i = 0,\ldots,N$, on the diagonal blocks. Again assuming the observation errors are distributed according to a normal distribution with error covariance matrix $R_i$ we define the cost function 
\beq
\label{eq:4dvarp}
J(x_0)=\frac{1}{2}\|y-\mathcal{H}(x)\|^2_{R^{-1}} +\frac{1}{2}\|x_0-x_0^B\|_{B^{-1}}^2=\frac{1}{2}\sum_{i=0}^N\|y_i-\mathcal{H}_i(x_i)\|^2_{R_i^{-1}} +\frac{1}{2}\|x_0-x_0^B\|_{B^{-1}}^2.
\eeq
In addition we introduce a discrete model for the evolution of the underlying  physical system from time $t_i$ to time $t_{i+1}$, described by the dynamical system equations
\[
x_{i+1} = \mathcal{M}_i(x_i),
\]
where $\mathcal{M}_i:\mathbb{R}^n\rightarrow\mathbb{R}^n$ can be time-dependent and nonlinear, describing the evolution of the dynamical system. Usually this forward operator requires the solution of a time-dependent partial differential equation. The problem
\begin{align}
\underset{x_0\in\mathbb{R}^n}{\operatorname{argmin }} J(x_0), \quad\text{where}\quad J(x_0)&= \frac{1}{2}\sum_{i=0}^N\|y_i-\mathcal{H}_i(x_i)\|^2_{R_i^{-1}} +\frac{1}{2}\|x_0-x_0^B\|_{B^{-1}}^2,\label{eq:4dvar} \\
\text{subject to the nonlinear model dynamics}\quad x_{i+1} &= \mathcal{M}_i(x_i) \label{eq:4dvar1}
\end{align}
is a constrained optimisation problem and is called \emph{4D-Var}, four dimensional variational data assimilation. Here, the subscripts refer to the time index. Note that in (\ref{eq:4dvar}) the regularisation term only involves the initial state $x_0$ at time $t_0$. The minimisation provides the initial condition of the model that most closely fits the data.

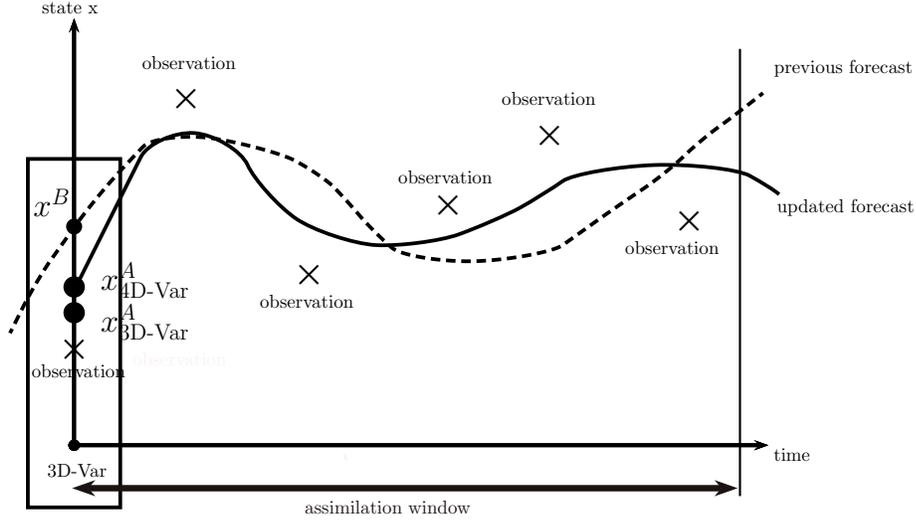
\begin{figure}[h!]
\begin{center}
\psscalebox{0.71 0.71} 
{
\begin{pspicture}(0,-6.937188)(19.050526,4.188)
\definecolor{colour0}{rgb}{0.98039216,0.9529412,0.9529412}
\definecolor{colour1}{rgb}{0.10980392,0.07058824,0.07058824}
\rput(2.9,-5.562812){\psaxes[linecolor=black, linewidth=0.07, tickstyle=bottom, axesstyle=axes, labels=none, ticks=none, dx=0.85cm, dy=1.0cm, showorigin=false]{*->}(0,0)(0,0)(13,8)}
\pscustom[linecolor=black, linewidth=0.07, linestyle=dashed, dash=0.17638889cm 0.10583334cm]
{
\newpath
\moveto(1.7,-3.462812)
\lineto(1.95,-2.962812)
\curveto(2.075,-2.712812)(2.325,-2.287812)(2.45,-2.112812)
\curveto(2.575,-1.937812)(2.825,-1.587812)(2.95,-1.412812)
\curveto(3.075,-1.237812)(3.375,-0.862812)(3.55,-0.662812)
\curveto(3.725,-0.46281198)(4.05,-0.11281198)(4.2,0.037188012)
\curveto(4.35,0.18718801)(5.0,0.26218802)(5.5,0.18718801)
\curveto(6.0,0.11218801)(6.75,-0.087811984)(7.0,-0.21281199)
\curveto(7.25,-0.33781198)(7.625,-0.587812)(7.75,-0.712812)
\curveto(7.875,-0.837812)(8.125,-1.112812)(8.25,-1.262812)
\curveto(8.375,-1.412812)(8.675,-1.712812)(8.85,-1.862812)
\curveto(9.025,-2.012812)(9.825,-2.137812)(10.45,-2.112812)
\curveto(11.075,-2.087812)(11.875,-1.937812)(12.05,-1.812812)
\curveto(12.225,-1.687812)(12.575,-1.437812)(12.75,-1.312812)
\curveto(12.925,-1.187812)(13.325,-0.912812)(13.55,-0.76281196)
\curveto(13.775,-0.612812)(14.15,-0.312812)(14.3,-0.16281198)
\curveto(14.45,-0.012811986)(14.775,0.26218802)(14.95,0.38718802)
\curveto(15.125,0.512188)(15.425,0.73718804)(15.8,1.037188)
}
\rput[bl](16.010527,1.337188){previous forecast}
\rput[bl](16.0,-5.862812){time}
\rput[bl](2.3,2.537188){state x}
\psdots[linecolor=black, dotstyle=x, dotsize=0.4](2.9,-3.762812)
\psdots[linecolor=black, dotstyle=x, dotsize=0.4](7.3,-2.362812)
\psdots[linecolor=black, dotstyle=x, dotsize=0.4](11.8,0.23718801)
\psdots[linecolor=black, dotstyle=x, dotsize=0.4](14.421053,-1.362812)
\psdots[linecolor=black, dotstyle=x, dotsize=0.4](5.0,0.937188)
\psdots[linecolor=black, dotstyle=x, dotsize=0.4](9.9,-1.062812)
\pscustom[linecolor=black, linewidth=0.07]
{
\newpath
\moveto(2.9,-2.662812)
\lineto(3.15,-2.162812)
\curveto(3.275,-1.912812)(3.525,-1.412812)(3.65,-1.162812)
\curveto(3.775,-0.912812)(4.025,-0.412812)(4.15,-0.16281198)
\curveto(4.275,0.08718801)(4.775,0.312188)(5.15,0.28718802)
\curveto(5.525,0.26218802)(6.025,-0.037811987)(6.15,-0.312812)
\curveto(6.275,-0.587812)(6.725,-1.087812)(7.05,-1.312812)
\curveto(7.375,-1.537812)(8.125,-1.787812)(8.55,-1.812812)
\curveto(8.975,-1.837812)(9.75,-1.737812)(10.1,-1.612812)
\curveto(10.45,-1.487812)(11.025,-1.237812)(11.25,-1.112812)
\curveto(11.475,-0.987812)(11.875,-0.737812)(12.05,-0.612812)
\curveto(12.225,-0.48781198)(13.075,-0.33781198)(13.75,-0.312812)
\curveto(14.425,-0.287812)(15.325,-0.387812)(15.55,-0.51281196)
\curveto(15.775,-0.63781196)(16.025,-0.787812)(16.1,-0.862812)
}
\rput[bl](16.063158,-1.2733383){updated forecast}
\rput[bl](10.910526,0.8161354){observation}
\rput[bl](8.973684,-0.6522857){observation}
\rput[bl](4.173684,1.4792932){observation}
\rput[bl](6.3789473,-2.9838645){observation}
\psframe[linecolor=black, linewidth=0.07, dimen=outer](3.8,-0.16281198)(2.0,-6.762812)
\pscustom[linecolor=colour0, linewidth=0.07]
{
\newpath
\moveto(8.0,-5.6628118)
}
\rput[bl](4.0,-4.062812){\textcolor{colour0}{observation}}
\rput[bl](2.1052632,-4.2838647){observation}
\rput[bl](2.4,-6.1628118){3D-Var}
\pscustom[linecolor=colour0, linewidth=0.12]
{
\newpath
\moveto(7.957895,-5.799654)
\lineto(8.0,-5.7786007)
}
\psline[linecolor=colour1, linewidth=0.12, arrowsize=0.05291667cm 2.1,arrowlength=1.43,arrowinset=0.3]{<->}(2.863158,-6.3680754)(15.326316,-6.3680754)
\rput[bl](7.2210526,-6.852286){\textcolor{colour1}{assimilation window}}
\rput[bl](13.221052,-1.9891278){observation}
\psline[linecolor=black, linewidth=0.04](15.368421,1.8635038)(15.368421,-6.515444)
\psdots[linecolor=black, dotsize=0.3](2.9,-1.462812)
\psdots[linecolor=black, dotsize=0.4](2.9,-3.080994)
\psdots[linecolor=black, dotsize=0.4](2.9,-2.595392)
\rput[bl](2.142106,-1.2733383){\LARGE{$x^B$}}
\rput[bl](3.4,-2.808266){\LARGE{$x^A_{\text{4D-Var}}$}}
\rput[bl](3.4,-3.590847){\LARGE{$x^A_{\text{3D-Var}}$}}
\end{pspicture}
}
\end{center}
\caption{Illustration of 3D-Var and 4D-Var for state estimation in one dimension.}
\label{fig:varpic}
\end{figure} 
Figure \ref{fig:varpic} illustrates the difference between 4D-Var and 3D-Var (for illustration purposes the state dimension is $n=1$). The 3D-Var problem finds the best estimate using an observation and a background state $x^B$ at one specific time point (here at the initial time), the best estimate is called $x^A_{\text{3D-Var}}$. In 4D-Var observations are obtained throughout a time window - often called the assimilation window - and the optimisation is done over a time window. The best estimate, taking into account the background state, that is the previous forecast and the observations within the assimilation window is computed, where we have to take the dynamical forecast model into account. The best estimate at the initial state is then $x^A_{\text{4D-Var}}$, which of course is not necessary the same as $x^A_{\text{3D-Var}}$. Both estimates can be used to evolve the model forward and to obtain a better (updated) forecast.

The cost function in (\ref{eq:3dvar}) or (\ref{eq:4dvar}) is minimised using an iterative method (e.g. steepest descent, conjugate gradient, a quasi-Newton method etc. \cite{Nocedal2006}). We choose $x_0^{(0)}$ (often $x_0^{(0)} = x_0^B$) and update 
\[
x_0^{(\ell+1)} = x_0^{(\ell)}+\alpha^{(\ell)} d^{(\ell)}, \quad\ell = 0,1,2,\ldots,
\]
where $\alpha^{(\ell)}$ is a damping parameter obtained through line search, for example, and $\ell$ is the iteration index. The search direction is $d^{(\ell)} = -\nabla J(x_0^{(\ell)})$ for gradient descent, or $d^{(\ell)} = -(\nabla^2 \tilde{J}(x_0^{(\ell)}))^{-1} \nabla J(x_0^{(\ell)})$ for a quasi-Newton method (where $\nabla^2 \tilde{J}(x_0^{(\ell)})$ is an approximation of the Hessian of $J(x_0^{(\ell)})$, or Jacobian matrix of $\nabla J(x_0^{(\ell)})$). 

Clearly, 4D-Var (\ref{eq:4dvar})-(\ref{eq:4dvar1}) is a constrained optimisation problem, and the gradient $\nabla J(x_0^{(\ell)})$ is obtained via an adjoint approach (see \cite{Asch2016}): We introduce Lagrange multipliers $\lambda_i\in\mathbb{R}^n$ at time $t_i$, $i=1,\ldots,N$ and define the Lagrangian 
\[
\mathcal{L}(x_i,\lambda_i) = J(x_0) +\sum_{i=0}^{N-1}\lambda_{i+1}^T(x_{i+1}-\mathcal{M}_i(x_i)).
\]
Necessary conditions for the minimum of the 4D-Var cost function subject to the constraint are then found by taking the variations of $\mathcal{L}$ with respect to $\lambda_i$ and $x_i$ (KKT conditions, see \cite{Nocedal2006})
\begin{align*}
\frac{\partial \mathcal{L}(x_i,\lambda_i)}{\partial x_0} &= B^{-1}(x_0-x_0^B) + H_0^TR_0^{-1}(\mathcal{H}_0(x_0)-y_0)-M_0^T\lambda_1 = 0,\\
\frac{\partial \mathcal{L}(x_i,\lambda_i)}{\partial x_i} &= H_i^TR_i^{-1}(\mathcal{H}_i(x_i)-y_i)-M_i^T\lambda_{i+1}+ \lambda_i  = 0,\quad i=1,\ldots,N,\\
\frac{\partial \mathcal{L}(x_i,\lambda_i)}{\partial \lambda_i} &= x_{i}-\mathcal{M}_{i-1}(x_{i-1}) = 0\quad i=1,\ldots,N,
\end{align*}
where $M_i\in\mathbb{R}^{n\times n}$ and $H_i\in\mathbb{R}^{p\times n}$ are the Jacobians of the forward operator $\mathcal{M}_i$ and the observation operator $\mathcal{H}_i$, evaluated at $x_i$, respectively, that is 
\[
M_i = \frac{\partial \mathcal{M}_i}{\partial x}(x_i),\quad H_i = \frac{\partial \mathcal{H}_i}{\partial x}(x_i).
\]
The adjoint equations for the adjoint variables $\lambda_i$, $i = 0,\ldots,N+1$, that measure the sensitivity of the cost function to changes in $x_i$, are then given by
\begin{align}
\lambda_{N+1}&=0\\
\lambda_i &= M_i^T\lambda_{i+1}-H_i^TR_i^{-1}(\mathcal{H}_i(x_i)-y_i),\quad i = N,\ldots,0.\label{eq:adjoint}
\end{align}
They provide an efficient method to compute the gradient of the objective function (\ref{eq:4dvar}), which is then given by
\[
\nabla J(x_0)=-\lambda_0+B^{-1}(x_0-x_0^B).
\]
The general method for solving the 4D-Var optimisation problem is sketched in Algorithm \ref{alg:4dvar}.
\begin{algorithm}
\caption{4D-Var}\label{alg:4dvar}
\begin{algorithmic}
  \State \textbf{Input:} error covariance matrices $R_i$ and $B$, routines to apply model and observation operators $\mathcal{M}_i$ and $\mathcal{H}_i$ and their linearisations $M_i$ and $H_i$, respectively, and observations $y_i$ for $i=0,\ldots,N$, maximum number of iterations 
$\ell_{\max}$.
  \State Initialise the iteration $\ell=0$ and $x_0^{(0)}=x_0^B$.
   \While{$\|\nabla J(x_0^{(\ell)})\|>\varepsilon$ or $\ell\le \ell_{\max}$}
   \State Compute cost function $J(x_0^{(\ell)})$ using the forward model.
   \State Compute $\nabla J(x_0^{(\ell)})$ using the adjoint equations.
   \State Apply descent method (gradient descent, Newton, quasi-Newton $\ldots$) to compute descent direction $d^{(\ell)}$.
   \State Update the initial condition $x_0^{(k+1)} = x_0^{(\ell)}+\alpha^{(\ell)} d^{(\ell)}$.
   \State Set $\ell=\ell+1$.
\EndWhile 
\end{algorithmic}
\end{algorithm}

The problem in (\ref{eq:4dvar})-(\ref{eq:4dvar1})  is often called \emph{strong constraint 4D-Var}. In the presence of uncertainty we can model the imperfect state dynamics by
\[
x_{i+1}=\mathcal{M}_i(x_i)+\eta_i,
\]
where $\eta_i\sim\mathcal{N}(0,Q_i)$ is the model error which is assumed to be Gaussian with error covariance $Q_i$, uncorrelated in time and uncorrelated with the background and observation errors. The relaxation of the strong constraint is commonly used in sequential data assimilation as we will see in Section \ref{sec:seq}. For variational data assimilation a weak constraint was proposed in \cite{Sasaki1970}, however, due to high computational costs is not used heavily in practice. In the past decades, with the availability of increasing computing power, there has been greater interest in this method \cite{Fisher2011a, Freitag2018,Gratton2018,Tremolet2007a}. The cost function for weak constraint 4D-Var is given by
\beq
\label{eq:4dvarw}
J(x)= \frac{1}{2}\sum_{i=1}^N\|y_i-\mathcal{H}_i(x_i)\|^2_{R_i^{-1}} +\frac{1}{2}\sum_{i=0}^{N-1}\|x_{i+1}-\mathcal{M}_i(x_{i})\|^2_{Q_{i+1}^{-1}}  +\frac{1}{2}\|x_0-x_0^B\|_{B^{-1}}^2,
\eeq
where $x = [x_0^T,\ldots,x_N^T]^T$, and $x_i$ is the model state at time step $t_i$, $i=0,\ldots,N$. The resulting optimisation problem 
\beq
\label{eq:4dvarw1}
\text{argmin}_{x\in\mathbb{R}^{n(N+1)}} J(x)
\eeq 
can become extremely large (in three spatial and the additional time dimension) and computationally expensive to solve, as $\mathcal{M}_i$ usually arises from a the discretisation of a nonlinear partial differential equation. We will discuss more detailed solution approaches to \eqref{eq:4dvarw1} in Section \ref{sec:lin}.

Note that the variational approach assumes that prior and likelihood have Gaussian distributions. As such the posterior is also Gaussian if the observation operators $h$ (for 3D-Var) or $\mathcal{H}$ (for 4D-Var, in which case it also includes the model operator) are linear. These are quite strong assumptions and, by minimising the cost function $J$ in (\ref{eq:3dvar}), (\ref{eq:4dvarp}), or (\ref{eq:4dvar})-(\ref{eq:4dvar1})  we only compute the posterior mean (the MAP estimator). 

We are able to approximate the posterior covariance matrix $A$ by computing the inverse of the Hessian, $A = (\nabla^2 J(x))^{-1}$. This covariance matrix is exact for linear models and Gaussian distributions, in which case it is given by $A = (B^{-1}+\mathrm{H}^T\mathrm{R}^{-1}\mathrm{H})^{-1}$, with $\mathrm{H}=H$ and $\mathrm{R}=R$ for 3D-Var, and $\mathrm{H}=[H_0,H_1M_0,H_2 M_1 M_0,\ldots,H_N M_{N-1}\cdots M_0]^T$ and $\mathrm{R} = \text{diag}(R_0,R_1,\ldots,R_N)$ for 4D-Var, respectively (see more details about this at the end of Section \ref{sec:seq} and in reference \cite{gejadze2018}). If the models are nonlinear or the distributions non-Gaussian, the inverse of the Hessian $(\nabla^2 J(x))^{-1}$ is only an approximation of the posterior covariance matrix.

To conclude this section we emphasise that variational data assimilation is essentially based on optimal control theory, that is we minimise a cost function subject to a constraint arising from the state dynamics. The minimisation requires techniques from numerical optimisation. For more insight into the relation of variational data assimilation and PDE-constrained optimisation we refer to \cite{Asch2016, Fisher2009,griffith1997}. 

\subsection{Sequential data assimilation and Kalman filters}
\label{sec:seq}

In sequential data assimilation we correct the model state estimate whenever observations are available, that is we incorporate the observations sequentially. Recall the Bayesian inference problem (\ref{eq:bayes}) with $x_i\in\mathbb{R}^n$ and $y_i\in\mathbb{R}^p$ observations at time $t_i$, for $i = 0,\ldots,N$. Given a prior uncertainty $\pi_X(x_i)$, find the updated conditional posterior uncertainty $\pi_{X|Y}(x_i|y_i)$ of $x_i$, given the observation $y_i$. Using the likelihood $\pi_{Y|X}(y_i|x_i)$ and Bayes' formula we would like to compute
\beq
\label{eq:bayes1}
\pi_{X|Y}(x_i|y_i)\propto\pi_{Y|X}(y_i|x_i)\pi_X(x_i).
\eeq
Again, ideally we wish to compute the the entire probability density function $\pi_{X|Y}(x_i|y_i)$ at every time step $t_i$.
We consider a discrete-time linear system dynamics,
\begin{align}
x_{i+1} &= M_i x_i+w_i,\label{eq:dyn}\\
y_{i} &= H_ix_i+v_i,\label{eq:obs}
\end{align}
where $w_i\sim\mathcal{N}(0,Q_i)$ and $v_i\sim\mathcal{N}(0,R_i)$ represent model and observation/measurement errors and are assumed to be independent, and Gaussian with error covariance matrices $Q_i\in\mathbb{R}^{n\times n}$ and $R_i\in\mathbb{R}^{p\times p}$ respectively. Hence the errors are uncorrelated in time and between each other. Moreover, model and observation operators are assumed to be linear here, however, extensions to nonlinear models exist.
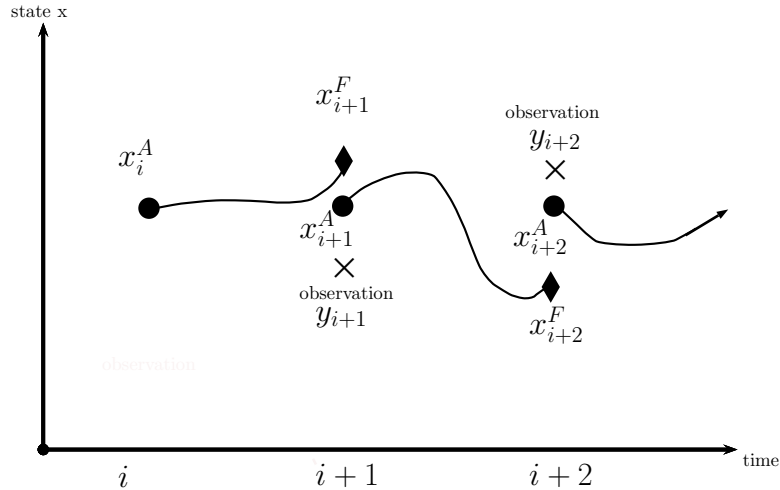
\begin{figure}[h!]
\begin{center}
\psscalebox{0.71 0.71} 
{
\begin{pspicture}(0,-6.1740303)(16.8,5.1740303)
\definecolor{colour0}{rgb}{0.98039216,0.9529412,0.9529412}
\rput(2.9,-5.8259697){\psaxes[linecolor=black, linewidth=0.07, tickstyle=bottom, axesstyle=axes, labels=none, ticks=none, dx=0.85cm, dy=1.0cm, showorigin=false]{*->}(0,0)(0,0)(13,8)}
\rput[bl](16.0,-6.12597){time}
\rput[bl](2.3,2.2740302){state x}
\psdots[linecolor=black, dotsize=0.4](4.8789473,-1.3049172)
\psdots[linecolor=black, dotstyle=x, dotsize=0.4](8.542106,-2.4154437)
\psdots[linecolor=black, dotstyle=x, dotsize=0.4](12.489473,-0.5891278)
\pscustom[linecolor=colour0, linewidth=0.07]
{
\newpath
\moveto(8.0,-5.92597)
}
\rput[bl](4.0,-4.3259697){\textcolor{colour0}{observation}}
\pscustom[linecolor=colour0, linewidth=0.12]
{
\newpath
\moveto(7.957895,-6.062812)
\lineto(8.0,-6.0417585)
}
\rput[bl](4.305263,-0.6733383){\LARGE{$x_i^A$}}
\rput[bl](4.305263,-6.52){\LARGE{$i$}}
\rput[bl](8.0,-6.52){\LARGE{$i+1$}}
\rput[bl](12.0,-6.52){\LARGE{$i+2$}}
\pscustom[linecolor=black, linewidth=0.04]
{
\newpath
\moveto(4.9894733,-1.304917)
\lineto(5.326316,-1.2417597)
\curveto(5.494737,-1.2101806)(5.9,-1.1680756)(6.1368423,-1.1575489)
\curveto(6.373684,-1.1470221)(6.889474,-1.1628119)(7.168421,-1.1891278)
\curveto(7.447368,-1.2154437)(7.831579,-1.1838641)(7.936842,-1.1259699)
\curveto(8.042105,-1.0680755)(8.215789,-0.94702274)(8.28421,-0.88386476)
\curveto(8.352632,-0.82070684)(8.452632,-0.6943909)(8.48421,-0.6312329)
\curveto(8.515789,-0.568075)(8.531579,-0.46807498)(8.484211,-0.35754883)
}
\psdots[linecolor=black, dotstyle=diamond*, dotsize=0.4](8.526316,-0.42070672)
\psdots[linecolor=black, dotstyle=diamond*, dotsize=0.4](12.4,-2.7786014)
\rput[bl](8,0.4635038){\LARGE{$x_{i+1}^F$}}
\rput[bl](7.7,-3){observation}
\rput[bl](8,-3.5364962){\LARGE{$y_{i+1}$}}
\rput[bl](11.557895,0.37929326){observation}
\rput[bl](12,-0.220706723){\LARGE{$y_{i+2}$}}
\psdots[linecolor=black, dotsize=0.4](8.505263,-1.262812)
\psdots[linecolor=black, dotsize=0.4](12.463158,-1.262812)
\pscustom[linecolor=black, linewidth=0.04]
{
\newpath
\moveto(8.505263,-1.2207068)
\lineto(8.810526,-0.9786017)
\curveto(8.963158,-0.85754883)(9.357894,-0.6943909)(9.6,-0.6522858)
\curveto(9.842106,-0.6101801)(10.157895,-0.6417597)(10.23158,-0.71544373)
\curveto(10.3052635,-0.7891278)(10.468422,-1.0259699)(10.557896,-1.1891278)
\curveto(10.647368,-1.3522857)(10.789473,-1.6680756)(10.842105,-1.8207068)
\curveto(10.894736,-1.973338)(11.005263,-2.2522857)(11.063158,-2.3786018)
\curveto(11.121052,-2.5049176)(11.268421,-2.689128)(11.357896,-2.7470222)
\curveto(11.447369,-2.804917)(11.584211,-2.8891277)(11.631579,-2.9154437)
\curveto(11.678947,-2.9417596)(11.8052635,-2.9786017)(11.884212,-2.9891279)
\curveto(11.96316,-2.999654)(12.073684,-2.9786017)(12.105264,-2.9470227)
\curveto(12.136844,-2.9154437)(12.205264,-2.8628118)(12.315789,-2.799654)
}
\pscustom[linecolor=black, linewidth=0.04]
{
\newpath
\moveto(12.589474,-1.3259699)
\lineto(12.747368,-1.473338)
\curveto(12.826316,-1.5470221)(12.957895,-1.673338)(13.010528,-1.7259699)
\curveto(13.063159,-1.7786018)(13.178947,-1.8733381)(13.2421055,-1.9154438)
\curveto(13.3052635,-1.9575489)(13.673684,-1.9943908)(13.978948,-1.9891279)
\curveto(14.284211,-1.9838648)(14.668421,-1.9364966)(14.747368,-1.8943915)
\curveto(14.826316,-1.8522857)(14.952632,-1.7838647)(15.0,-1.7575488)
\curveto(15.047368,-1.731233)(15.1,-1.699654)(15.115789,-1.6838646)
}
\psline[linecolor=black, linewidth=0.05, arrowsize=0.05291667cm 2.0,arrowlength=1.4,arrowinset=0.0]{->}(14.926316,-1.8101804)(15.7263155,-1.3259699)
\rput[bl](12,-3.8417594){\LARGE{$x_{i+2}^F$}}
\rput[bl](7.7,-2.038646){\LARGE{$x_{i+1}^A$}}
\rput[bl](11.7,-2.2022856){\LARGE{$x_{i+2}^A$}}
\end{pspicture}
}
\end{center}
\caption{Illustration of sequential data assimilation for the Kalman Filter in one dimension.}
\label{fig:kfpic}
\end{figure}
The Kalman filter (as well as many other sequential data assimilation schemes) then follow two steps: A forecast or prediction step and an analysis or correction step, illustrated in Figure \ref{fig:kfpic}. At time step $t_i$ we have an analysis $x_i^A$ available, which we assume to be the best estimate of the state at time $t_i$ (we have incorporated all our knowledge about observations and previous forecasts at time $t_i$).
Then the state dynamics are evolved forward in time using (\ref{eq:dyn}) producing a forecast $x_{i+1}^F$, which becomes the background state at step $i+1$. We also have a set of observations $y_{i+1}$ available which we combine with $x_{i+1}^F$ in order to get the best estimate of the state, a so-called analysis, $x_{i+1}^A$. This analysis is obtained using the \text{best linear unbiased estimate} (BLUE), see \cite{humpherys2012}. This is again used to evolve the state forward using the model dynamics. 

Let $x_i^t\in\mathbb{R}^n$ be the unknown exact true state of the system. Then we define the error in the forecast at time $t_i$ by $e_i^F = x_i^F-x_i^t$, and similar, the error in the analysis at time $t_i$ by $e_i^A = x_i^A-x_i^t$. Moreover, the forecast and analysis error covariance matrices can then be expressed as 
\[
P_i^F = \text{cov}(e_i^F) = \mathbb{E}[e_i^F (e_i^F)^T]\quad\text{and}\quad P_i^A = \text{cov}(e_i^A) = \mathbb{E}[e_i^A (e_i^A)^T],
\]
respectively.
\paragraph{Forecast/predictor step} Given a previous analysis state estimate $x_i^A$ at time $t_i$, an estimate of $x_{i+1}^t$ at time $t_{i+1}$ is given by the application of the model dynamics (\ref{eq:dyn})
\[
x_{i+1}^F = M_i x_i^A.
\] 
For the corresponding error covariance matrix $P_{i+1}^F=\mathbb{E}[e_{i+1}^F (e_{i+1}^F)^T]$ we observe $e_{i+1}^F = M_i e_i^A + w_i$. Since the model is linear and the analysis error $e_i^A$ and model error $w_i$ are assumed to be uncorrelated we obtain (using properties of the expected value and $\mathbb{E}[w_iw_i^T] = Q_i$)
\begin{equation}
\label{eq:propcov}
P_{i+1}^F = M_i P_i^A M_i^T + Q_i.
\end{equation}
\paragraph{Analysis/corrector step} 
Given an a-priori estimate $x_i^F$ and observations $y_i$ at time step $i$, the goal of the analysis step is to compute the optimal a-posteriori estimate, $x_i^A$ as a linear combination of $x_i^F$ and $y_i$ of the form
\begin{equation}
\label{eq:kalmangain}
x_i^A = x_i^F + K_i(y_i-H_ix_i^F),
\end{equation}
where $K_i\in\mathbb{R}^{n\times p}$ is called the Kalman gain matrix \cite{Kalman1960} and the general form of (\ref{eq:kalmangain}) arises from the assumption that the estimate $x_i^A$ of the true vector $x_i^t$ should be both linear and unbiased. The a-posteriori error covariance is given by $P_i^A = \mathbb{E}[e_i^A (e_i^A)^T] = \mathbb{E}[(K_i(-H_i e_i^F +v_i) + e_i^F)(K_i(-H_i e_i^F +v_i) + e_i^F)^T]$,
where we have used the observation equation (\ref{eq:obs}). With the definition of $P_i^F:=\mathbb{E}[e_i^F (e_i^F)^T]$ and $\mathbb{E}[v_iv_i^T] =: R_i$, we obtain 
\beq
\label{eq:pia}
P_i^A = (I-K_iH_i)P_i^F(I-K_iH_i)^T+K_iR_iK_i^T,
\eeq
where we have also used that the forecast error $e_i^F$ and observation error $v_i$ are uncorrelated, random variables, that is $\mathbb{E}[e_i^F v_i^T]=0$.

The Kalman gain matrix $K_i$ in \eqref{eq:pia} is chosen to minimise the a-posterior variance $\text{tr}(P_i^A)$. We evaluate $\frac{\partial{P_i^A}}{\partial K_i}=0$, and using results from matrix differential calculus we obtain the Kalman gain 
\beq
\label{eq:kgain}
K_i = P_i^F H_i^T (H_i P_i^F H_i^T + R_i)^{-1},
\eeq
for details of the derivation, see \cite{Asch2016,humpherys2012,bardsley2019}. Substituting $K_i$ into \eqref{eq:pia} then yields
\beq
\label{eq:covA}
P_i^A = (I-K_iH_i)P_i^F.
\eeq
Note that, for $P_i^F$ small (or zero) the Kalman gain is also small (or zero) and the a-posteriori estimate $x_i^A$ is heavily geared towards $x_i^F$. Similarly, if the observation error covariance matrix $R_i$ is zero (for perfect observations) the estimate $x_i^A$ is steered towards the observations (in particular, for $H_i$ square and nonsingular $x_i^A=H_i^{-1}y_i$). We remark that the error covariance matrices in the Kalman filter equations satisfy discrete algebraic Riccati equations \cite{lancaster1995,gustafsson2011some}.

Clearly the Kalman Filter has some shortcomings, as it is only optimal for linear models and observation operators, and Gaussian error statistics (in which case the mean and error covariances, which are propagated in the Kalman Filter, are sufficient to describe the probability density function of the state estimates).

Several methods have been proposed to overcome these issues, some we mention here, for more details we refer to \cite{Asch2016}. The \textit{Extended Kalman Filter} (EKF) 
applies to systems of the form
\begin{align*}
x_{i+1} &= \mathcal{M}_i (x_i)+w_i,\\
y_{i} &= \mathcal{H}_i(x_i)+v_i,
\end{align*}
where $\mathcal{M}_i$ and $\mathcal{H}_i$ are nonlinear. The nonlinearities are dealt with in the filter equations by linearising the model and observation operator about $x_i^F$ and $x_i^A$, in the computations of the Kalman gain and error covariance matrices, that is $M_i = \frac{\partial\mathcal{M}_i}{\partial x}(x_i^A)$, $H_i = \frac{\partial\mathcal{H}_i}{\partial x}(x_i^F)$. For computing the analysis state estimate and the forecast state estimate the nonlinear operators are used. We obtain an approximation to the best linear unbiased estimate, see \cite{Gelb1974,Jazwinski1970}. An algorithmic description of the EKF is given in Algorithm \ref{alg:EKF}.
\begin{algorithm}
\caption{Extended Kalman filter}\label{alg:EKF}
\begin{algorithmic}
  \State \textbf{Input:} error covariance matrices $R_i$ and $Q_i$, routines to apply model and observation operators $\mathcal{M}_i$ and $\mathcal{H}_i$ and their linearisations $M_i$ and $H_i$, respectively, and observations $y_i$ for $i=0,\ldots,N$.
  \State Initialise the system state $x_0^F$ and the corresponding error covariance matrix $P_0^F$.
   \For{$i=0,\ldots,N$}
   \State Compute Kalman gain $K_i = P_i^F H_i^T (H_i P_i^F H_i^T + R_i)^{-1}$.
   \State Compute state estimate $x_i^A = x_i^F + K_i(y_i-\mathcal{H}_i(x_i^F))$.
   \State Compute error covariance estimate $P_i^A = (I-K_iH_i)P_i^F$.
   \State Compute the forecast state $x_{i+1}^F = \mathcal{M}_i(x_i^A)$.
   \State Compute the forecast error covariance $P_{i+1}^F = M_i P_i^A M_i^T + Q_i$.
\EndFor
\end{algorithmic}
\end{algorithm}

There are a lot of similarities between the Kalman filter discussed in this section and variational data assimilation discussed in Section \ref{sec:var}. In particular, one can interpret the analysis/corrector step in the Kalman filter in a variational form: the Kalman filter state  estimate $x_i^A$ from \eqref{eq:kalmangain} with the Kalman gain \eqref{eq:kgain} exactly represent the equations we obtained for computing the solution to the 3D-Var problem in Section \ref{sec:var}. Hence, for the variational formulation of the Kalman filter analysis step we have 
\begin{equation}
\label{eq:varform}
x_i^A = \text{argmin}_{x_i\in\mathbb{R}^n} J_i(x_i), \quad\text{where}\quad J_i(x_i) =  \left(\frac{1}{2}\|y_i-H_i(x_i)\|^2_{R_i^{-1}} +\frac{1}{2}\|x_i-x_i^F\|_{(P_i^F)^{-1}}^2\right),
\end{equation}
which is a quadratic problem for a linear observation operator $H_i$ and hence has a unique solution which can be obtained by setting the gradient $\nabla J_i(x_i)=0$. For a nonlinear observation operator the solution may not be unique. The Hessian of this cost function (for linear observation operator $H_i$) is given by
\[
\nabla^2 J_i(x_i) = H_i^T R_i^{-1}H_i + (P_i^F)^{-1},
\] 
and, with the Sherman-Morrison-Woodbury formula we have 
\[
\nabla^2 J_i(x_i)^{-1} = P_i^F-P_i^F H_i^T (H_i P_i^FH_i^T+R_i )^{-1} H_iP_i^F,
\]
which is precisely the a-posteriori error covariance matrix $P_i^A$ given in \eqref{eq:covA}. Therefore, the inverse of the Hessian is equal to the posterior covariance (in the linear case with Gaussian errors). For nonlinear problems, $\nabla^2 J_i(x_i)^{-1}$ is an approximation to the posterior covariance. Solving \eqref{eq:varform} iteratively is more advantageous than computing the Kalman gain directly for very large problems with sparsity structure. More details on this idea can be found in \cite{bardsley2019}.

To complete this section we add a remark about Kalman smoothers. Filtering algorithms make use of observations as they become available and provide the best estimate of a state given all past information and the current observation. Smoothing algorithms provide the best estimate of a system, using past, present and future information. So the Kalman smoother is equivalent to the Kalman filter at the final time step. Moreover, it can be shown that for linear problems with linear observation operator and model dynamics, and the same initial background error covariance, the Kalman filter and 4D-Var result in the same state estimate for the same time step (when the same observations have been used), see  \cite{FrePot2013,Li2001,Fisher2005}.

In the next sections we give an overview of important approaches and contributions from numerical linear algebra to solve variational and sequential data assimilation problems efficiently and provide an outlook for challenges ahead.

\section{Solutions to the optimisation problem arising in variational data assimilation}
\label{sec:lin}

Variational data assimilation leads to large PDE-constrained optimisation problems. We will first concentrate on 4D-Var \eqref{eq:4dvar}-\eqref{eq:4dvar1} and later discuss weak constraint 4D-Var \eqref{eq:4dvarw}.

\paragraph{Incremental 4D-Var and Gauss-Newton method}

Since the full nonlinear minimisation \eqref{eq:4dvar}-\eqref{eq:4dvar1} is difficult to solve, special algorithms for optimisation and linear algebra are required. In geoscience applications an incremental approach \cite{Courtier1994} was proposed, which is merely a Gauss-Newton method for nonlinear least squares problems \cite{Lawless2005,Lawless2005inv,Gratton2007}. In incremental variational data assimilation the solution to the nonlinear optimisation problem is approximated by a sequence of minimisations of quadratic (and hence convex) cost functions, which are obtained by linearising both the model and the observation operators. Let $x_0^{(\ell)}$ be the $\ell$th estimate to the solution at $x_0$. We linearise the cost function \eqref{eq:4dvar} around the model trajectory from this estimate: we linearise the model and observation operators about $x_0^{(\ell)}$ and then obtain the next iterate by the increment $\delta x_0^{(\ell)}$,
\begin{equation}
\label{eq:newtonup}
x_0^{(\ell+1)} =x_0^{(\ell)}+\delta x_0^{(\ell)}.
\end{equation}
We use this expansion and substitute it into the nonlinear cost function \eqref{eq:4dvar}, which we then linearise about the model trajectory obtained from $x_0^{(\ell)}$. We then see that the increment $\delta x_0^{(\ell)}$ is obtained by minimising the quadratic problem (the incremental cost function)
\begin{equation}
\label{eq:quadcost}
\tilde{J}^{(\ell)}(\delta x_0^{(\ell)}) = \frac{1}{2} \sum_{i = 0}^{N} \|H_i \delta x_i^{(\ell)} - d_i^{(\ell)} \|^2_{R_i^{-1}} + \frac{1}{2}\|\delta x_0^{(\ell)} - (x_0^B-x_0^{(\ell)})\|^2_{B^{-1}}, 
\end{equation}
where $d_i^{(\ell)} = y_i - \mathcal{H}_i(x_i^{(\ell)})$ and $x_i^{(\ell)}$ is the nonlinear trajectory computed from the current estimate at the initial time $x_0^{(\ell)}$, using the nonlinear model trajectory. $H_i$ is the linearisation of the observation operator $\mathcal{H}_i$ about $x_i^{(\ell)}$. The perturbations $\delta x_i^{(\ell)}$ satisfy the linearised constraint
\[
\delta x_{i+1}^{(\ell)} = M_i(\delta x_i^{(\ell)}),
\]
where $M_i$ is the linear solution operator of the nonlinear model $\mathcal{M}_i$ linearised around the nonlinear trajectory. We therefore obtain a so-called inner-outer iterative method. The outer iteration (with iteration index $\ell$) is represented by the update of the nonlinear trajectory \eqref{eq:newtonup}. The minimisation of the quadratic cost function \eqref{eq:quadcost} is the inner iteration, it essentially amounts to the solution of a large linear system. At each iteration the forward model, and its linearisation, are evaluated in order to compute the cost function \eqref{eq:quadcost}, and the adjoint model \eqref{eq:adjoint} is applied in order to compute the gradient of the cost function. 

Incremental variational data assimilation can be shown to be a version of the Gauss-Newton method applied to the original nonlinear cost function \eqref{eq:4dvarp}. In order to illustrate this, consider a general nonlinear least squares problem
\beq
\label{eq:nonl}
\min_{x}\Phi(x) = \frac{1}{2}f(x)^T f(x) = \frac{1}{2}\|f(x)\|^2,
\eeq
with $f:\mathbb{R}^n\rightarrow \mathbb{R}^q$ a twice continuously differentiable function $[f_1(x),\ldots,f_q(x)]$ and $\|\cdot\|$ denoting the Euclidean norm. Let $G(x) = f'(x)$ be the $q\times n$ Jacobian of $f(x)$. The gradient and Hessian of $\Phi(x)$ are then given by
\[
\nabla\Phi(x) = G(x)^T f(x),\quad \nabla^2\Phi(x) = G(x)^TG(x)+\sum_{j=1}^q f_j(x)\nabla^2 f_j(x).
\]
The Gauss-Newton method discards the expensive second term in the Hessian when applying Newton's method to $\nabla\Phi(x) = G(x)^T f(x)=0$. 
\begin{algorithm}
\caption{Gauss-Newton method}\label{alg:gn}
\begin{algorithmic}
\State \textbf{Input:} Routines to compute $f(x)$ and its Jacobian $G(x)$, maximum number of iterations $\ell_{\max}$.
\State Initialise the iteration $\ell=0$ and $x^{(0)}=x_0$
\While{$\|G(x^{(\ell)})^Tf(x^{(\ell)})\|>\varepsilon$ or $\ell\le \ell_{\max}$}
   \State Solve $G(x^{(\ell)})^TG(x^{(\ell)}) \delta x^{(\ell)}=-G(x^{(\ell)})^Tf(x^{(\ell)})$.
   \State Update $x^{(\ell+1)}= x^{(\ell)}+\delta x^{(\ell)}$.
   \State Set $\ell=\ell+1$.
\EndWhile 
\end{algorithmic}
\end{algorithm}
The Gauss-Newton method is described in Algorithm \ref{alg:gn}. Note that the second step of the Algorithm is equivalent to solving the linearised least squares problem 
\beq
\label{eq:nonlingn}
\min_{s\in\mathbb{R}^n}\|G(x^{(\ell)}) s + f(x^{(\ell)})\|
\eeq
at each iteration $\ell$. If we define 
\[
f(x_0) = \mat{c}B^{-1/2}(x_0-x_0^B)\\ R_0^{-1/2}(y_0-\mathcal{H}_0(x_0))\\\vdots\\R_N^{-1/2}(y_N-\mathcal{H}_N(x_N))\rix
\]
subject to $x_{i+1} = \mathcal{M}_i(x_i)$, then the general cost function \eqref{eq:nonl}, with $x=x_0$ is equivalent to the 4D-Var cost function \eqref{eq:4dvar}-\eqref{eq:4dvar1}, where here $q=n+(N+1)p$. When applying the Gauss-Newton method to \eqref{eq:nonl}, then the linearised least squares problem \eqref{eq:nonlingn} is equivalent to finding the solution to the quadratic cost function \eqref{eq:quadcost} with $x=x_0$. More details and also suitable stopping criteria for the inner-outer iteration arising within the Gauss-Newton method were discussed in \cite{Lawless2005,Lawless2005inv,Lawless2006,Gratton2007}. 

Besides Gauss-Newton, there are of course standard methods to solve nonlinear least squares problems, such as Levenberg-Marquardt and Quasi-Newton methods (in particular BFGS), see, for example \cite{Fisher2009,Nocedal2006}. 

\paragraph{The inner iteration and preconditioning}

The quadratic cost function \eqref{eq:quadcost} (see also \eqref{eq:nonlingn}) within the Gauss-Newton (incremental) approach can be minimised using a conjugate gradient (CG) method \cite{hestenes1952,Golub2012,
Nocedal2006,gill2019}. As the problem is usually very large, in practice, the inner loop, the CG iteration, is applied to a system with lower spatial resolution and with a preconditioner \cite{Fisher2009}. One of the earliest approaches is the multilevel setting in \cite{Courtier1994}, where the quadratic cost function is replaced by a lower resolution model, in order to obtain a multi-resolution scheme. 

In addition, two-level preconditioning is usually applied within 4D-Var in the following way. We consider the cost function
\begin{equation}
\label{eq:34dvar}
J(x_0)=\frac{1}{2}\|y-\mathcal{H}(x)\|^2_{R^{-1}} +\frac{1}{2}\|x_0-x_0^B\|_{B^{-1}}^2,
\end{equation}
where $\mathcal{H}:[x^T_0,\ldots, x^T_N]^T\rightarrow [\mathcal{H}_0(x_0)^T,\ldots, \mathcal{H}_N(x_N)^T]^T$ for 4D-Var and 
$\mathcal{H}:x_0\rightarrow \mathcal{H}_0(x_0)$, where $\mathcal{H}_0(x_0)=h(x_0)$ for 3D-Var, respectively. The formulation in \eqref{eq:34dvar} aims to treat 3D-Var and 4D-Var simultaneously. 

At each step of the Gauss-Newton process an equation of the form $G(x^{(\ell)})^TG(x^{(\ell)}) \delta x^{(\ell)}=-G(x^{(\ell)})^Tf(x^{(\ell)})$ needs to be solved, which results in
\[
(B^{-1}+H^T R^{-1} H) \delta x^{(\ell)} = B^{-1}(x_0^B-x^{(\ell)})+H^T R^{-1} (y-\mathcal{H}(x^{(\ell)})),
\]
where $H$ depends on $x^{(\ell)}$. Note that $H$ is the linearisation of $\mathcal{H}$ at $x^{(\ell)}$. For 4D-Var $H$ includes the model operator, $H = [H_0,H_1 M_0,\ldots,H_N M_{N-1}\cdots M_0]^T\in\mathbb{R}^{(N+1)p\times n}$ and $R = \mathrm{diag}(R_0, R_1, \cdots, R_N)\in \mathbb{R}^{(N+1)p \times (N+1)p}$, for 3D-Var $H = H_0$ and $R=R_0$. We exclude the superscript $\ell$ for ease of notation. The Hessian of the 4D-Var cost function at the $\ell$th linearised problem is given by $B^{-1}+H^TR^{-1}H$ \cite{gejadze2008}. 

The first level preconditioning employs a linear change in variables,
\[\delta x = L \delta \tilde x,\quad\text{where}\quad B = LL^T,\]
that is, the Cholesky factor of the background error covariance matrix. The Hessian of the cost function then becomes $I+L^TH^TR^{-1}HL$, which ensures that the smallest eigenvalue of the transformed Hessian is one. Moreover, since the rank of $L^TH^TR^{-1}HL$ is smaller than the dimension of the system, there are many unit eigenvalues. Note that another option is to employ a transformation with $B$ to obtained the transformed Hessian  $I+BH^TR^{-1}H$, which is no longer symmetric (but the use of this transformation may be necessary if the factorisation of $B$ is not available). Both transformed Hessians share the same eigenvalues \cite{el2013}; they are all greater than or equal to one. Detailed analysis on the conditioning of the Hessian can be found in \cite{haben2011,tabeart2018}.

At a second preconditioning level the spectrum of the (symmetric) Hessian is used. This utilises the fact that a few of the dominant eigenvalues and corresponding eigenvectors can be obtained using the Lanczos method \cite{Lanczos1950,cullum2002}, and, one can compute those Hessian eigenpairs within the CG iteration itself \cite{Golub2012}. After $k$ steps of the CG algorithm, a few approximate leading eigenpairs $(\lambda_i,w_i)$, $i=1,\ldots,k$ of the Hessian are available and one can approximate the Hessian by
\beq
\label{eq:approxH}
C = I+\sum_{i=1}^k (\lambda_i-1)w_i w_i^T,
\eeq
where $k$ is much smaller than the dimension of the linear system. Hence, after one step of the Gauss-Newton method (the outer iteration), approximations to the Hessian, $C$, are available for the next step of the Gauss-Newton iteration, and preconditioners $C^{-1}$, or $C^{-\frac{1}{2}}$ for $H$ can be obtained by replacing $\lambda_i$ in \eqref{eq:approxH} by $1/\lambda_i$ or $1/\sqrt{\lambda_i}$, respectively, see  \cite{anderson2000,morales2000,Tremolet2007}. In \cite{el2013} this procedure was extended to non-symmetric Hessians by using the bi-conjugate gradient method.  

The work in \cite{ramage2016} takes up the idea of approximating the Hessian by eigenpairs obtained from the Lanczos procedure. Matrix-vector products with the Hessian are expensive (since they require the evaluation of the linearised forward and adjoint models), so even obtaining the limited memory representation \eqref{eq:approxH} may be computationally infeasible for large $k$. Therefore \cite{ramage2016} propose a multilevel version of the limited memory Hessian, where the eigenvalue decompositions are obtained from several coarser levels and fed through to the finer level, which enhances the algorithm in terms of reducing the number of matrix vector products at the fine grid levels. In addition, multigrid solvers and multigrid preconditioners for the solution of the variational data assimilation problem were considered recently in  \cite{Debreu2015,gratton2013b}.

The limited memory preconditioners (LMP) discussed above are often referred to as spectral LMP, as they use the spectral information of the Hessian approximation \eqref{eq:approxH}. More general versions of LMP were investigated in \cite{tshimanga2008} and \cite{Gratton2011}. In the latter the authors also show the equivalence of a certain reduced version of 4D-Var and the SEEK filter, discussed in Section \ref{sec:kfred}, and use this equivalence to accelerate the convergence of the Gauss-Newton method.

\paragraph{The dual formulation of 4D-Var}

As observed in the previous paragraph, both for 3D-Var and for 4D-Var, each inner iteration of the Gauss-Newton method essentially attempts to minimise the cost function 
\beq
\label{eq:primal}
\tilde{J}(\delta x) =\frac{1}{2}\delta x^T B^{-1}\delta x + \frac{1}{2}(H\delta x-d)^TR^{-1} (H\delta x-d),
\eeq
see \eqref{eq:quadcost}, which is equivalent to solving 
\[
(B^{-1}+H^T R^{-1} H) \delta x = H^T R^{-1}d\quad\text{or}\quad  \delta x = (B^{-1}+H^T R^{-1} H)^{-1} H^T R^{-1}d.
\]
Here we have neglected all sub- and superscripts for simplicity. 
The minimisation of the cost function in \eqref{eq:primal} is often referred to as the primal approach, as the minimisation takes place in model space. Using the Sherman-Morrison-Woodbury formula we write the solution as 
\[
 \delta x  = BH^T (HBH^T+R)^{-1}d,
\]
which can be obtained by solving the smaller system in observation space $(HBH^T+R)\lambda = d$ and setting $\delta x  = BH^T\lambda$. This method is known as dual formulation of 3D-Var/4D-Var, or PSAS (Physical-space Statistical Analysis System) \cite{courtier1997}. It is easy to see that the cost function for the dual formulation is given by
\[
D(\lambda) =\frac{1}{2}\lambda ^T (HBH^T+R)\lambda -\lambda^T d.
\]
If the dimension of the observation space $p$ is significantly smaller than that of the model state space $n$, the dual formulation can reduce both memory usage and computational cost compared to the primal approach. Preconditioned conjugate gradient methods for the dual approach were considered in \cite{gratton2009,gratton2013} and convergence properties of the primal and dual approaches were investigated in \cite{el2008,akkraoui2010}. Most recently, so called $B$-preconditioned minimisation algorithms for variational data assimilation were introduced in the paper \cite{gurol2014b}, which also contains a good literature review on the dual formulation and PSAS.

\paragraph{Weak constraint 4D-Var and saddle point formulation of the inner iteration}

Weak constraint 4D-Var requires minimisation over all state variables within the assimilation window and is therefore more computationally expensive. The incremental approach \cite{Courtier1994} for the more general weak constraint 4D-Var cost function \eqref{eq:4dvarw} can be formulated as follows. We approximate the 4D-Var cost function by a quadratic function of an increment
\begin{equation}
\delta x^{(\ell)} = x^{(\ell+1)} - x^{(\ell)},
\end{equation}
where ${x^{(\ell)}} = \left[(x_0^{(\ell)})^T, (x_1^{(\ell)})^T, \cdots, (x_N^{(\ell)})^T\right]^T$ denotes the $\ell$th iterate of the Gauss-Newton algorithm applied to weak-constraint 4D-Var. This increment $\delta x^{(\ell)}$ is a solution to the minimisation of the linearised cost function
\begin{equation}
\label{eq:WeakIncP}
\begin{aligned}
\tilde{J}^{\ell}(\delta x^{(\ell)}) &= \frac{1}{2} \| \delta x_0^{(\ell)} - b_0^{(\ell)} \|^2_{B^{-1}} + \frac{1}{2} \sum_{i= 0}^{N} \| d_i^{(\ell)} - H_i \delta x_i^{(\ell)} \|^2_{R_i^{-1}} + \frac{1}{2} \sum_{i = 0}^{N-1} \|\delta x_{i+1}^{(\ell)} - M_{i} \delta x_{i}^{(\ell)} -c_{i+1}^{(\ell)}\|^2_{Q_{i+1}^{-1}},
\end{aligned}
\end{equation}
where $M_i$ and $H_i$, are linearisations of $\mathcal{M}_i$ and $\mathcal{H}_i$ about the current state trajectory $x^{(\ell)}$, and  
$b_0^{(\ell)} = x_0^b - x_0^{(\ell)}$, $d_i^{(\ell)} = y_i - \mathcal{H}_i(x_i^{(\ell)})$, and $c_{i+1}^{(\ell)} = \mathcal{M}_i (x_{i}^{(\ell)}) - x_{i+1}^{(\ell)}$. Dropping the superscript for the $\ell$th iterate for simplicity, the linearised cost function \eqref{eq:WeakIncP} can be written more concisely as 
\begin{equation}
\label{eq:WeakInccondensed}
\tilde{J}(\delta x) = \frac{1}{2} \| \mathbf{L} \delta x - b \|^2_{\mathbf{D}^{-1}} + \frac{1}{2} \|d - \mathbf{H} \delta x \|^2_{\mathbf{R}^{-1}},
\end{equation}
where $\delta x = \left[\delta x_0^T, \delta x_1^T, \cdots, \delta x_N^T\right]^T$ and $\mathbf{L}$ and $\mathbf{H}$ are matrices of size ${(N+1)n \times (N+1)n}$ and ${(N+1)p \times (N+1)n}$, respectively:
\begin{equation}
\label{eq:LH}
\mathbf{L} = \begin{bmatrix}
I 	 &  	  & 	   & \\
-M_0 & I  	  & 	   & \\
& \ddots & \ddots & \\
&		  & -M_{N-1}   & I
\end{bmatrix}, \quad \mathbf{H} = \begin{bmatrix}
H_0 &&& \\ &H_1&& \\ &&\ddots& \\ &&&H_N,
\end{bmatrix}
\end{equation}
which can be thought of as all-at-once model and observation operators over the assimilation window. Here we assume $y_i\in\mathbb{R}^p$, but this can be generalised to $y_i\in\mathbb{R}^{p_i}$. We assume there is no correlation between the errors at each time steps, and hence the covariance matrices are block diagonal matrices
\beq
\mathbf{D} = \mathrm{diag}(B,Q_1, \cdots, Q_N)\in \mathbb{R}^{(N+1)n \times (N+1)n},\quad\text{and}\quad \mathbf{R}= \mathrm{diag}(R_0, R_1, \cdots, R_N)\in \mathbb{R}^{(N+1)p \times (N+1)p}.
\eeq
Moreover, the vectors $b$ and $d$ are given by 
\begin{align*} 
b = \left[b_0^T, c_1^T, \cdots, c_N^T\right]^T \in \mathbb{R}^{(N+1)n},\quad\text{and}\quad d = \left[d_0^T, d_1^T, \cdots, d_N^T\right]^T \in \mathbb{R}^{(N+1)p}.
\end{align*}
The system above can be written as a saddle point problem \cite{Fisher2011,Fisher2011a,Fisher2016,Fisher2017}, a form that recently has seen a lot of interest for data assimilation problems, see also \cite{Freitag2018,green2019}. With new variables 
\[\lambda = \mathbf{D}^{-1}(b-\mathbf{L}\delta x)\quad\text{and}\quad \mu = \mathbf{R}^{-1}(d-\mathbf{H}\delta x),\]
the gradient of the cost function \eqref{eq:WeakInccondensed} provides a constraint and altogether the coupled linear system
			\begin{equation}
\label{eq:saddle}
			\begin{bmatrix}
			\mathbf{D} & 0 & \mathbf{L} \\
			0 & \mathbf{R} & \mathbf{H} \\
			\mathbf{L}^T & \mathbf{H}^T & 0
			\end{bmatrix} \begin{bmatrix}
			\lambda \\ \mu \\ \delta x
			\end{bmatrix} = \begin{bmatrix}
			b \\ d \\ 0
			\end{bmatrix},
			\end{equation}
a very large, sparse, symmetric indefinite saddle point system needs to be solved at every inner iteration. A vast amount of literature on saddle point problems and their solution via Krylov methods and preconditioners is available, see for example \cite{Benzi2005} and references therein. For this particular saddle point problem low-rank limited memory preconditioners exploiting the structure of the saddle point problem were proposed and analysed in \cite{Fisher2016} (see also \cite{gratton2018note}). In the work \cite{Freitag2018,green2019} the Kronecker structure of the saddle point problem was used in order to compute low-rank solutions to GMRES. The convergence of the saddle point formulation of weak constrained 4D-Var was reviewed in \cite{Gratton2018}, and spectral estimates for the saddle point system were obtained in \cite{dauickait2019}.

\paragraph{Hybrid methods and inexact approaches}

It has been observed that incremental weak constraint 4D-Var for minimising the large scale cost function \eqref{eq:4dvarw} may diverge, hence  the authors in \cite{mandel2016hybrid} add a regularisation term and thereby replace the Gauss-Newton approach by the Levenberg-Marquardt method. In addition they use an ensemble Kalman smoother (see Section \ref{sec:kfred} for a discussion on ensemble methods) within the minimisation process and thereby apply a hybrid approach. Such hybrid methods combining variational and sequential ensemble approaches are popular as ensemble methods are naturally parallelisable and do not require adjoint operators \cite{clayton2013operational}. 

A parallel-in-time approach for solving the strong constraint 4D-Var optimisation problem was proposed in \cite{rao2016time}. 

When both the model operator and the gradient are not available exactly, inexact methods need to be used. In \cite{bergou2016,bellavia2018} a Levenberg-Marquardt method is proposed and investigated for dealing with inexact gradients and Jacobians.

\section{The Kalman filter and low-rank approximations}
\label{sec:kfred}

The Kalman filter is impractical for large dimensional systems. It requires the storage and evolution of large covariance matrices $P_i^{A/F}$, which are both not feasible for very high dimensional systems (for example, in oceanography and numerical weather prediction, the state dimension is of the order of $10^7$ and higher). The propagation of the error covariance in \eqref{eq:propcov} requires a number of integrations of the forward model equal to the dimension of the system. Moreover matrix inversion within the Kalman gain computation \eqref{eq:kgain} is expensive. Hence a range of approximate Kalman filters have been developed for large systems, either by using a simplified or reduced order model \cite{dee1991,Cohn1996,Farrell2001} to propagate the covariance matrices (that is, to propagate the error statistics) or by using a reduced state space or error space \cite{cane1996mapping,Pham1998,Verlaan1997}. Many of the approaches are quite similar and most of them rely on low-rank approximations of the error covariance matrices. We discuss some methods below. 

\paragraph{Reduced-rank Kalman filters}

The singular evolutive extended (SEEK) Kalman filter algorithm (see, for example \cite{Pham1998,Brasseur2006,Rozier2007,livings2008}) is one of the best known reduced rank square root (RRSQRT) filters. It is assumed that the covariance matrices $P$ arising within the algorithm have low-rank form and can be written as $P = SS^T$, where $S\in\mathbb{R}^{n\times r}$ with $r\ll n$. The factorisation can be obtained via a truncated eigenvalue decomposition, for example. The Kalman filter equations are then rewritten using the matrices $S_i^F$ and $S_i^A$, the low rank approximations of the forecast error and analysis error covariance matrix, respectively. The equation for the Kalman gain \eqref{eq:kgain} becomes
\[
K_i = S_i^F (H_i S_i^F)^T (H_iS_i^F (H_iS_i^F)^T + R_i)^{-1},
\]
or, using the Sherman-Morrison Woodbury identity 
\[K_i = S_i^F [I_r + (H_i S_i^F)^T R_i^{-1} H_iS_i^F]^{-1} (H_iS_i^F)^TR_i^{-1}.\]
Note that often $R_i$ is a (block) diagonal matrix, and the latter version of the low rank Kalman gain can be computed at lower cost in $r\ll n$ dimensions. The analysis increment is $x_i^A - x_i^F = K_i(y_i-H_i(x_i^F))$, and therefore a linear combination of the columns of $S_i^F$. Substituting the low rank Kalman gain into \eqref{eq:covA} yields $P_i^A = S_i^A (S_i^A)^T$ with 
\[
S_i^A = S_i^F [I_r + (H_i S_i^F)^T R_i^{-1} H_iS_i^F]^{-1/2},
\] 
where the inverse of the square root is taken in the lower dimensional space of dimension $r$. 

For the forecast step, the propagation of the error covariance matrix is done via 
\[
P_{i+1}^F = \tilde{S}_{i+1}^F (\tilde{S}_{i+1}^F)^T + Q_i, \quad\text{where}\quad \tilde{S}_{i+1}^F = M_i S_i^A,
\]
or, for nonlinear models, via the finite difference approximation
\beq
\label{eq:sif}
\{\tilde{S}_{i+1}^{F}\}_{\ell} = \mathcal{M}_i(x_i^A+ \{S_i^A\}_{\ell})-\mathcal{M}_i(x_i^A),\quad\ell = 1,\ldots,r,
\eeq
where $\{\cdot\}_{\ell}$ refers to the $\ell$th column. In order to write $P_{i+1}^F = S_{i+1}^{F} (S_{i+1}^F)^T$, and conserving the rank $r$ some assumptions need to be made about $Q_i$ (see \cite{Verron1999,Brasseur2006} for details), otherwise a rank reduction, for example via computing an SVD, may be required at every step in order to keep the rank of the covariance matrices small.

\begin{algorithm}
\caption{SEEK filter}\label{alg:seek}
\begin{algorithmic}
  \State \textbf{Input:} error covariance matrices $R_i$ and $Q_i$, routines to apply model and observation operators $\mathcal{M}_i$ and $\mathcal{H}_i$ and their linearisations $M_i$ and $H_i$, respectively, and observations $y_i$ for $i=0,\ldots,N$.
  \State Initialise the system state $x_0^F$ and the corresponding error covariance matrix in low-rank form $P_0^F = S_0^F (S_0^F)^T$.
   \For{$i=0,\ldots,N$}
   \State Compute low rank Kalman gain $K_i = S_i^F [I_r + (H_i S_i^F)^T R_i^{-1} H_iS_i^F]^{-1} (H_iS_i^F)^TR_i^{-1}$.
   \State Compute state estimate $x_i^A = x_i^F + K_i(y_i-\mathcal{H}_i(x_i^F))$.
   \State Compute error covariance estimate $P_i^A = S_i^A (S_i^A)^T$, where $S_i^A =S_i^F [I_r + (H_i S_i^F)^T R_i^{-1} H_iS_i^F]^{-1/2}$.
   \State Compute the forecast state $x_{i+1}^F = \mathcal{M}_i (x_i^A)$.
   \State Compute the forecast error covariance $P_{i+1}^F = \tilde{S}_{i+1}^F (\tilde{S}_{i+1}^F)^T+ Q_i$ where $\tilde{S}_{i+1}^F$ is given by \eqref{eq:sif}.
\EndFor
\end{algorithmic}
\end{algorithm}

The SEEK filter, shown in Algorithm \ref{alg:seek} is just one example of a reduced rank square root filter (RRSQRT), see \cite{Verlaan1997,Chandrasekar08}. Another example, shown in Algorithm \ref{alg:RRSQRT} is a more general version of the SEEK filter and makes sure the rank of the error covariance matrix does not increase during the iteration. Again, we drop the time index, and assume a low rank approximation of the error covariance matrix  $P = SS^T$ is possible, where $S\in\mathbb{R}^{n\times r}$ with $r\ll n$. The Kalman gain can then be written as 
\[
K_i = S_i^F (L_i^F)^T (L_i^F (L_i^F)^T + R_i)^{-1},\quad\text{where}\quad L_i^F = H_i S_i^F, 
\]
and simple computations (again using the Sherman Morrison Woodbury identity), show that 
\[
P_i^A = S_i^F \left(I_r -(L_i^F)^T(L_i^F (L_i^F)^T+R_i)^{-1} L_i^F\right)(S_i^F)^T,
\]
and using a matrix square root (of a smaller $r\times r$ matrix), one can write $S_i^A =  S_i^F (I_r-(L_i^F)^T(L_i^F (L_i^F)^T+R_i)^{-1} L_i^F)^{\frac{1}{2}}$. Several algorithms are available for computing the matrix square root, e.g. \cite{higham1997,druskin1998}, but often a Cholesky factorisation is used. After the analysis step the dimension of the system is reduced, keeping only $r-s$ eigenmodes of $(S_i^A)^T S_i^A$, where $s<r$. This avoids an increase of the rank of the covariance matrix when variability through the model error (with covariance matrix $Q_i$) is introduced. The details of the RRSQRT filter are given in Algorithm \ref{alg:RRSQRT}. 
\begin{algorithm}
\caption{Reduced rank square root filter (RRSQRT)}\label{alg:RRSQRT}
\begin{algorithmic}
  \State \textbf{Input:} error covariance matrices $R_i$ and $Q_i$, routines to apply model and observation operators $\mathcal{M}_i$ and $\mathcal{H}_i$ and their linearisations $M_i$ and $H_i$, respectively, and observations $y_i$ for $i=0,\ldots,N$.
  \State Initialise the system state $x_0^F$ and the corresponding error covariance matrix in low-rank form $P_0^F = S_0^F (S_0^F)^T$ and decompose $Q_i = T_iT_i^T$, a rank $s$ factorisation of $Q_i$, with $s<r$.
   \For{$i=0,\ldots,N$}
   \State Compute low rank Kalman gain $K_i = S_i^F (L_i^F)^T ({L_i^F}(L_i^F)^T+R_i)^{-1}$, where $L_i^F=H_iS_i^F$.
   \State Compute state estimate $x_i^A = x_i^F + K_i(y_i-\mathcal{H}_i(x_i^F))$.
   \State Compute the low rank factor of error covariance $S_i^A =  S_i^F (I_r-(L_i^F)^T(L_i^F (L_i^F)^T+R_i)^{-1} L_i^F)^{\frac{1}{2}}$. 
\State Compute an eigenvalue decomposition (or low-rank factorisation) $V\Lambda V^T =(S_i^A)^T S_i^A$.
   \State Select largest $r-s$ eigenvalues and corresponding eigenvectors $\tilde{V}:=V(:,1:r-s)$, set $\tilde{S}_i^A = S_i^A\tilde{V}$.
   \State Compute the forecast state $x_{i+1}^F = \mathcal{M}_i (x_i^A)$.
   \State Compute low rank factor $S_{i+1}^F = [M_i\tilde{S}_i^A, T_i]$ of the forecast error covariance $P_{i+1}^F = S_{i+1}^F (S_{i+1}^F)^T$.
\EndFor
\end{algorithmic}
\end{algorithm}
Reduced rank filters are cheaper to implement than the full filtering algorithm, by using low rank approximations of covariance matrices. Only $r$ forward model integrations are necessary. However, they still use linearisations of the nonlinear operators $\mathcal{M}_i$ and $\mathcal{H}_i$ in order to propagate the error covariance matrices. The ensemble Kalman filter overcomes this issue and utilises the full nonlinear model to propagate the covariances.

\paragraph{Ensemble Kalman filter}

The need to propagate a probability distribution is an important feature of the ensemble Kalman Filter (EnKF). It is also a major challenge as the propagation of large covariance matrices of size $n$ is very expense. The ensemble Kalman filter is also a reduced rank method, as it requires the propagation and analysis of a small number of ensemble members. It was proposed in \cite{Evensen1994} (see also \cite{HoutekamerMitchell98,Evensen2009,Bishop2001,Tippett2003}) and is essentially a Monte Carlo implementation of the Bayesian update. There are plenty of variants of the ensemble Kalman filter, with the same idea behind all of them, the difference is in the implementation detail. We describe two versions briefly, the stochastic ensemble Kalman filter (EnKF) and the ensemble transform Kalman filter (ETKF), where, in the latter case, the linear algebra is performed in the ensemble subspace which is of much smaller dimension than the state space or the observation space. 

\begin{algorithm}
\caption{Ensemble Kalman Filter (EnKF)}\label{alg:enkf}
\begin{algorithmic}
  \State \textbf{Input:} error covariance matrices $R_i$, routines to apply forward model and observation operators $\mathcal{M}_i$ and $\mathcal{H}_i$, respectively, and observations $y_i$ for $i=0,\ldots,N$.
  \State Initialise the ensemble states $x_{k,0}^F$ where $k=1,\ldots,r$ (via random perturbations from the initial conditions $x_0^F$, for example).
   \For{$i=0,\ldots,N$}
   \State Perturb observations $y_{k,i} = y_i+v_k$, where $v_k\sim\mathcal{N}(0,R_i)$, $k=1,\ldots,r$.
   \State Compute the ensemble means
\[
\bar{x}_i^F = \frac{1}{r}\sum_{k=1}^r x_{k,i}^F,\quad
\bar{v} = \frac{1}{r}\sum_{k=1}^r v_k,\quad
\bar{y}_i^F = \frac{1}{r}\sum_{k=1}^r \mathcal{H}_i(x_{k,i}^F) 
\]
\State Compute the rectangular normalised ensemble matrices
\[
[X^F]_{k,i} = \frac{x_{k,i}^F-\bar{x}_i^F}{\sqrt{r-1}},\quad [L^F]_{k,i} = \frac{\mathcal{H}_i(x_{k,i}^F)-\bar{y}_i^F-(v_k-\bar{v})}{\sqrt{r-1}}.
\]
\State Compute the (approximate Kalman) gain: $K_i = X_i^F (L_i^F)^T (L_i^F(L_i^F)^T)^{-1}$.
\State Update the analysis ensemble
\[
x_{k,i}^A = x_{k,i}^F+K_i(y_{k,i}-\mathcal{H}_i(x_{k,i}^F)),\quad k=1,\ldots r.
\]
\State Compute the forecast ensemble $x_{k,i+1}^F = \mathcal{M}_i (x_{k,i}^A)$, $\quad k=1,\ldots r$.
\EndFor
\end{algorithmic}
\end{algorithm}

Suppose we have $r$ ensemble members, or prior samples, $\{x_k^F\}_{k=1}^r$. Note that here the subscript $k$ denotes the ensemble index, we neglect the time index in this explanation for simplicity. The forecast error covariance can be estimated using the empirical covariance 
\[
P^F = \frac{1}{r-1}\sum_{k=1}^r (x_k^F-\bar{x}^F)(x_k^F-\bar{x}^F)^T, \quad\text{where}\quad \bar{x}^F = \frac{1}{r}\sum_{k=1}^r x_k^F,
\]
or
\begin{equation}
\label{label:lrP}
P^F = X^F (X^F)^T\quad\text{where}\quad [X^F]_k = \frac{x_k^F-\bar{x}^F}{\sqrt{r-1}},
\end{equation}
and $[X^F]_k$ the $k$th column of the $n\times r$ matrix $X^F$. Each of the ensemble members is then updated using \eqref{eq:kalmangain}: $x_k^A = x_k^F + K(y-\mathcal{H} (x_k^F))$ to obtain a posterior ensemble $[X^A]_k$, where  
\[
[X^A]_k = \frac{x_k^A-\bar{x}^A}{\sqrt{r-1}},\quad\text{with}\quad \bar{x}^A = \frac{1}{r}\sum_{k=1}^r x_k^A.
\]
When computing the sample posterior covariance $P^A =  X^A (X^A)^T$ it turns out this is underestimated compared to the BLUE in \eqref{eq:pia}. A way around this is to perturb the observation vector, $y_k = y + v_k $, where $v_k\sim\mathcal{N}(0,R)$ and set 
\[
[L^F]_k = \frac{Hx_k^F-H\bar{x}^F-(v_k-\bar{v})}{\sqrt{r-1}}.
\]
Then it can be shown \cite{Asch2016} that $[X^A]_k =[X^F]_k- K[L^F]_k$ results in the correctly estimated $P^A =  X^A (X^A)^T$. The Kalman gain can be computed using $K = X^F (L^F)^T (L^F(L^F)^T)^{-1}$. Moreover it can be shown that within the algorithm we only require the application of the nonlinear operator $\mathcal{H}$ to the ensemble members, rather than its linearised version $H$.

The full version of the (stochastic) EnKF is given in Algorithm \ref{alg:enkf}. Note that $k$ is the ensemble index, $i$ is, as before, the time index. For more details on the algorithm we refer to \cite{Asch2016,harlim2018}. An error analysis for the ensemble Kalman filter analysis step was performed in \cite{kovalenko2011}.

\begin{algorithm}
\caption{Ensemble Transform Kalman Filter (ETKF)}\label{alg:etkf}
\begin{algorithmic}
  \State \textbf{Input:} error covariance matrices $R_i$, routines to apply forward model and observation operators $\mathcal{M}_i$ and $\mathcal{H}_i$, respectively, and observations $y_i$ for $i=0,\ldots,N$. $V\in\mathbb{R}^{r\times r}$ an orthogonal matrix such that $V\mathbf{1}= \mathbf{1}$.
  \State Initialise the ensemble states $\{x_{k}^F\}_{k=1,\ldots,r}$ at time $i=0$ and set $E^F =[x_1^F,\ldots, x_r^F]$.
   \For{$i=0,\ldots,N$}
   \State Compute ensemble mean $\bar{x}^F = E^F\mathbf{1}/r$ and ensemble matrix $X^F = (E^F-\bar{x}^F\mathbf{1}^T)/\sqrt{r-1}$.
   \State Compute the observation mean  $\bar{y} = Y\mathbf{1}/r$ where $Y = \mathcal{H}_i(E^F)$.
   \State Compute normalised observation ensemble $\tilde{L}^F = R_i^{-\frac{1}{2}} (Y-\bar{y}\mathbf{1}^T)/\sqrt{r-1}$. 
\State Compute normalised innovation vector $d = R_i^{-\frac{1}{2}}(y_i-\bar{y})$ and set $W = (I_r+ (\tilde{L}^F)^T \tilde{L}^F)^{-1}$.
\State Compute ensemble space coefficient vector $w^A = W(\tilde{L}^F)^T d$.
\State Update state estimate ensemble $E^A = \bar{x}^F\mathbf{1}^T + X^F(w^A\mathbf{1}^T +\sqrt{r-1} W^{\frac{1}{2}}V)$.
\State Compute the forecast ensemble $E^F = \mathcal{M}_i (E^A)$.
\EndFor
\end{algorithmic}
\end{algorithm}

For the deterministic version of the ensemble Kalman filter \cite{Bishop2001}, the ensemble transform Kalman filter (ETKF) which is illustrated in Algorithm \ref{alg:etkf}, we again write the forecast error covariance in low rank form \eqref{label:lrP}, with ensembles $\{x_k^F\}_{k=1}^r$. It is then assumed that the state estimate $x^A$ is of the form $x^A = \bar{x}^F + X^F w^A$, where $w^A$ is a vector of coefficients in the small dimensional ensemble subspace $\mathbb{R}^r$, and $X^F\in\mathbb{R}^{n\times r}$. Using \eqref{eq:kalmangain}, the mean of the analysis vector is given by $x^A = \bar{x}^F + K(y-\mathcal{H}(\bar{x}^F))$. Hence, we obtain
\[
\bar{x}^F + X^F w^A = \bar{x}^F + K(y-\mathcal{H}(\bar{x}^F)),
\]
and with the low rank approximation $P^F = X^F (X^F)^T$ within the Kalman gain $K$, and using the Sherman-Morrison-Woodbury formula again, we derive 
\[
w^A = (I_r+ (L^F)^T R^{-1}L^F)^{-1} (L^F)^T R^{-1}(y-\mathcal{H}(\bar{x}^F)),\quad\text{where}\quad L^F =\mathcal{H}X^F,
\]
and hence the Kalman gain is computed in the low dimensional ensemble space. For computing the posterior ensemble covariance matrix we proceed as in the RRSQRT derivation. This yields 
\[
X^A = X^F (I_r+ (L^F)^T R^{-1}L^F)^{-\frac{1}{2}}V,
\]
where $V$ is an arbitrary orthogonal matrix. To ensure that $X^A \mathbf{1}=0$, that is, the updated perturbations are centred at $x^A$ (similar to $X^F \mathbf{1}=0$) it is sufficient for $V \mathbf{1}=\mathbf{1}$ to hold. Here $\mathbf{1}$ is the vector of all ones. The posterior ensemble is then given by
\[
x_k^A = x^A+\sqrt{r-1} X^F  \left[W^{\frac{1}{2}}V\right]_k = \bar{x}^F+X^F\left(w^A+\sqrt{r-1}\left[W^{\frac{1}{2}}V\right]_k\right),
\] 
where $W = (I_r+ (L^F)^T R^{-1}L^F)^{-1}$. Note that within the ETKF, all matrix inversions and matrix square roots are carried out in the lower dimensional ensemble subspace of dimension $r$. Moreover, with a small number of ensemble members, only $r$ applications of the expensive forward model $\mathcal{M}_i$ and the operator $\mathcal{H}_i$ are necessary.

Using RRSQRT filters and EnKFs results in the increments being confined in a subspace spanned by the columns of the low rank matrix . There are localisation methods which overcome these issues, for details we refer to \cite{Asch2016,ott2004,hunt2007} and references therein. A whole range of variations of the ensemble and square root filters have been developed over the years, see the references in \cite{Asch2016,nerger2005comparison}.

\paragraph{Iterative solvers within the Kalman filter}

As mentioned at the end of Section \ref{sec:basic}, a variational form of the Kalman filter can be formulated and state estimate and posterior covariance are given by minimum and inverse Hessian of a quadratic cost function, respectively \cite{bardsley2019}. When the preconditioned CG method is used, this leads to a conjugate gradient ensemble Kalman filter, which has been discussed in \cite{bardsley2013ens}. Low-rank approximations to covariance and inverse covariance matrices can be obtained by exploiting the connection between conjugate gradient and Lanczos iterations \cite{bardsley2013krylov}. Moreover, in a similar way, the use of limited memory BFGS \cite{Nocedal2006} within the Kalman filter estimate has been investigated in \cite{auvinen2009}.

\section{Model reduction and dimension reduction approaches}
\label{sec:modred}

Data assimilation problems are often large, in particular when the model is described by discretised time-dependent partial differential equations and when large amounts of data need to be assimilated. This leads to either a large optimisation problem for variational data assimilation, or, the solution to large linear systems, eventually, for the optimisation problem and the Kalman filter. Over the years several reduction techniques have been proposed for data assimilation problems (or, inverse problems) which we will describe here briefly.

We will collect results from two different ideas. Usually the state space dimension $n$ in data assimilation is very large. This leads to very large covariance matrices which need to be stored, manipulated and inverted, which is expensive. Therefore, one approach is the reduction of the dimension of these covariance matrices, usually by using low-rank approximations.

A second approach considers the expensive, PDE-based model operator $\mathcal{M}$ (or $\mathcal{M}_i$, $i=0,\ldots,N$). In order to solve the optimisation problem in 4D-Var, or, in order to apply sampling based methods such as the ensemble Kalman filter (or, more general Markov Chain Monte Carlo methods), this model has to be evaluated many times at several points in space, which is expensive. Therefore reducing the dimension of that model operator is a second way of applying dimension reduction.

\paragraph{Reduced rank methods - reducing the dimension of the covariance matrix}

A key idea in dimension reduction for Bayesian inference is to exploit the fact that in updating the prior to the posterior, some directions in the high-dimensional state space are more important than others. In the case of Gaussian posteriors, this fact is quantified by the rate of decay of the eigenvalues of the Hessian of the data misfit \cite{Flathetal2011,bui2012extreme}. The decay leads to a low-rank approximation of the Hessian (and hence the posterior covariance). Quantitative error analysis of approximation methods for the posterior distribution in Bayesian inference have been performed in \cite{Spantinietal2015}.

We have already explained several reduced rank approaches applied to the Kalman Filter in Section \ref{sec:kfred}, which are based on reduced rank covariance matrices, so we will not repeat details  here. 

In recent years, a whole range of hybrid methods were developed, which combine ensemble (and hence low-rank) Kalman filters with variational data assimilation. This can be done in several ways and we point to \cite{Bannister2017} for a review of those ideas. One such method (often referred to EnVar) uses the low rank covariance matrix that arises within the ensemble Kalman filter (EnKF), and applies it within variational data assimilation as the background error covariance matrix, $B = SS^T$, where $S\in\mathbb{R}^{n\times r}$ and $r\ll n$. The resulting system has the same form as preconditioned 4D-Var, however with a low-rank version of the Hessian, and hence the optimisation problem is solved in the reduced system dimension, see \cite{Gratton2011}.

As discussed above, it is known that for linear inverse problems, the inverse of the Hessian is an approximation of the posterior error covariance. From a control theoretic view point, there is an equivalence between the Hessian and the observability Gramian. Hessian based model reduction using this viewpoint was investigated in \cite{lieberman2013hessian,bashir2007hessian,bashir2008hessian}
 
In \cite{Benner2018}, the Kronecker product structure of a discretised linear partial differential operator and a low-rank Arnoldi method were used to approximate posterior error covariance matrices of Gaussian posteriors using the idea of low-rank Hessians \cite{Flathetal2011,bui2012extreme,Spantinietal2015}. 

\paragraph{Model order reduction applied to the forward model operator}

A second approach for reducing the cost within sequential and variational data assimilation is reducing the cost of the forward and adjoint models by applying reduced order models (ROMs). The aim of model order reduction (MOR) is to find models that approximate and reflect the dynamics of the underlying large-scale system accurately, in ways that enable the reduction process to be implemented efficiently.
There is a large number of MOR techniques available, and many of them have been very popular in the system theoretic community \cite{benner2017}. 

The earliest reduction approaches for 4D-Var suggest using a simplified operator or a coarser grid within the minimisation of incremental 4D-Var \cite{Ide1997,Tremolet2004}. 

In \cite{Farrell2001} balanced truncation \cite{Moore1981}, a control theoretic approach for MOR, is  applied to the Kalman filter. The linearised model and observation operators $M_i\in\mathbb{R}^{n\times n}$ and $H_i\in\mathbb{R}^{p\times n}$ are projected onto lower dimensional subspaces, 
\beq
\label{eq:proj}
\hat{M}_i = U^T M_i V\in\mathbb{R}^{r\times r},\quad \hat{H}_i = H_i V\in\mathbb{R}^{p\times r},
\eeq
where $U\in\mathbb{R}^{n\times r}$ and $V\in\mathbb{R}^{n\times r}$ are projection matrices satisfying $U^TV = I$, obtained through balanced truncation. Balanced truncation is a method that retains the dominant observable and reachable states, which are the important ones for the system dynamics (after transforming both state and observation equation so that reachable and observable states can be expressed in the same basis). In \cite{Farrell2001} it is assumed that the time-dependent system underlying the problem has a time-invariant dominant part on which balanced truncation is performed. 

MOR via balanced truncation was also proposed for incremental 4D-Var in \cite{Boess2008,Boess2011,Lawless2008,Lawless2008a}. Model and observation operators are projected as in \eqref{eq:proj}, $\delta x$ is restricted to $\delta \hat{x} = U^T\delta x\in\mathbb{R}^r$, $r\ll n$, and the background error covariance matrix, $B$, is projected onto $U^T B U$. The approach was extended to weak constraint 4D-Var in \cite{Freitag2019}.

In \cite{Durbiano2001,Robert2005} it is assumed that the initial state $x_0$ in 4D-Var is contained in a space of reduced dimension $r\ll n$ about the background state,
\[
x_0 = x_0^B+\sum_{i=1}^r c_i w_i,\quad r\ll n,
\]
where $c_i$ are real coefficients and $w_i$ are linearly independent vectors containing the variability in the system. Minimisation of the reduced cost function then takes place in the reduced space of dimension $r$. The authors in \cite{arcucci2019optimal} use a truncated SVD approach to solve the 4D-Var problem in reduced spaces. 

Proper orthogonal decomposition (POD) methods were applied to variational data assimilation by several authors \cite{vermeulen2006,Cao2007,Daescu2007,dimitriu2009}. Snapshots are taken at various time steps from the model trajectory. An SVD is then taken of the matrix of snapshots $X = [x_1,\ldots, x_r]$, that is, $X=U\Sigma V^T$. Finally, $x_0$ is then projected onto the POD space spanned by the $r$ most important left singular vectors, that is, the ones corresponding to the largest singular values. The number of those left singular vectors is chosen significantly smaller than the dimension of the state space. POD was applied in \cite{altaf2013} to the adjoint in order to compute a reduced model. The work in \cite{Stefanescu2015} refines this work on POD by considering the discrete empirical interpolation method (DEIM) within POD for nonlinear dynamical systems.

A challenge for applying MOR techniques within data assimilation is the nonlinearity and time-dependence of the forward model operator. One idea is to use a data driven and online approach applied to the ROMs, where updates to the basis vectors of the ROMs are computed from combinations of reduced solutions, snapshots and adjoint information \cite{peherstorfer2015dynamic}.  

Dimension reduction is an important tool for Bayesian inverse problems as one often has to perform Markov chain Monte Carlo (MCMC) sampling to access the posterior distribution. However, each MCMC sample requires an expensive forward model solve. In \cite{Cui2014,Cui2015,Cui2016}, methods for  dimension reduction for nonlinear Bayesian inverse problems were described. The goal of these method is to approximate the likelihood using a reduced model. In \cite{Elmanetal2019}, POD-DEIM is successfully applied to the likelihood function in order to reduce the cost of each MCMC draw. A review of multi-fidelity methods was recently published \cite{peherstorfer2018survey}. For additional ideas for dimension reduction methods within data assimilation we refer to \cite{Asch2016} and references therein.

\section{Bayesian inference and Tikhonov regularisation, and other aspects within data assimilation}
\label{sec:baytic}

One aspect we have not considered in detail in this article is the link between data assimilation and Tikhonov regularisation. We refer to \cite{FrePot2013} and references therein. A key task in the Bayesian inverse problem framework is finding an informative prior. This corresponds to finding a computationally efficient penalty term in the regularisation approach, see \cite{calv2018,calvetti2018bayes}. Most of the work in this area considers Gaussian priors and linear, but very large scale static inverse problems. In reference \cite{chung2017generalized} methods based on Golub-Kahan bidiagonalisation for computing solutions based on Tikhonov regularisation are used, which avoid computing inverses and square roots of large covariance matrices. The idea was generalised to dynamic inverse problems, which also allows the use of a wide class of spatio-temporal priors \cite{Chung_2018,saibaba2018quantifying}.

When designing an observation or sensor network for data assimilation another important aspect not considered here is the placement of sensors. An important tool for this is computing the so-called observation impact, which provides a measure for important or redundant information. Mathematically this results in analysing the sensitivity of states with respect to the data, and eventually in a large linear system to solve. For more information and solution methods, including low-rank approaches, we refer to \cite{tremolet2008computation,cioaca2013,cioaca2014} and references therein.

Data assimilation using different regularisation terms was considered in \cite{budd2011,rao2017robust}, efficient solution techniques for such approaches require ideas from numerical linear algebra. 

There is a whole range of ideas for data assimilation we have not considered here. For literature on particle filters for Bayesian inference, for example, we refer to \cite{Asch2016,chen2003bayesian} and references therein.

\section{Software}

Software development for data assimilation algorithms is very much driven by geoscientists. Several packages are available to download and test algorithms, these also often provide examples. We list a few here.

PDAF, a software environment (written in FORTRAN) for ensemble data assimilation was developed by researchers from the Alfred Wegener Institute for Polar and Marine Research and is freely available \cite{Nerger2013}. DATes is a recently developed data assimilation testing suite written in PYTHON \cite{attia2017}. Finally, DART \cite{dart2009} is a FORTRAN software package developed  and  maintained  at the  National  Center for Atmospheric  Research.

MATLAB code, which is applied to a range of simple examples, both for Kalman filters and variational data assimilation, is available in the book by Law et al \cite{LawStuZyg2015}. Further MATLAB codes are provided which accompany the book \cite{ReiCot2013}.

\section{Conclusions}

In the era of ``big data'' it is important to be able to process and evaluate data, gain insight from data, extract knowledge from data and make predictions from data. However, making important decisions based only on data can be misleading as they might be incomplete or erroneous. Therefore, it is often crucial to combine data with mathematical models. This approach leads to the important area of \emph{data assimilation} which is evolving rapidly and continuously. 

Data assimilation uses tools from many different fields of mathematics, such as statistics, optimisation, machine learning, numerical linear algebra, mathematical modelling and scientific computing, to name a few.

For example, classical optimisation tools in variational data assimilation are constantly adapted to faster algorithms and new supercomputers. Variational data assimilation typically only provides a point estimate, but no uncertainty quantification. However, approximations of 
a-posteriori covariance matrices can be computed using efficient numerical linear algebra techniques. 

On the other hand, statistical learning techniques, such as Kalman filters and ensemble methods are regularly seeing improvements. Those techniques do provide uncertainty quantification as they are merely Monte Carlo implementation of the Bayesian update.

We have seen that in the field of data assimilation, traditional computational scientific modelling meets the area of statistics and data science in order to produce new and better algorithms and methods.

In this review we have stated and explained the most established data assimilation methods, starting from the point of view of Bayesian inference. We have focused on problems arising within the numerical linear algebra for these methods, that is, the solution to linear systems, preconditioning, and matrix methods for filtering and model reduction, and thereby given a extended (but by no means complete) review of the existing literature in numerical methods for data assimilation.

With the increasing size and complexity of datasets and larger models after discretisation of partial differential equations, many data assimilation problems involve operations on matrices with millions or billions of elements. This amount of large matrices brings new computational challenges to classical numerical linear algebra algorithms, for example solving large systems of linear equations, large linear (and nonlinear) regression problems, constructing low-rank matrix approximation etc.

The  efficiency  of  these  linear  algebra  applications  is  essential for the performance of data assimilation methods.

Hence, this review shows that efficient numerical linear algebra techniques and matrix computation tools are crucial within the subject of data driven modelling and data assimilation, even more so when the data and models become even larger.

\paragraph{Acknowledgements} The research of the author has been partially funded by Deutsche Forschungsgemeinschaft (DFG) - SFB1294/1 - 318763901 (associated member).

\end{document}